\documentclass[3p]{elsarticle}

\usepackage{url} 
\usepackage{algorithm2e}
\usepackage{caption}
\usepackage{subcaption}
\usepackage{amsmath}
\usepackage{amssymb}
\usepackage{mdframed}
\usepackage{lineno}
\usepackage{soul} 
\usepackage{natbib}
\usepackage{cleveref}
\usepackage{xcolor}

\definecolor{darkgreen}{rgb}{0,0.7,0}

\newcommand{\Hmatrix}{$\mathcal{H}$-matrix}

\newcommand{\QIRvec}{\m{Q}_{\operatorname{IR}}}
\newcommand{\QIR}{Q_{\operatorname{IR}}}
\newcommand{\Qabsvec}{\m{Q}_{\operatorname{abs}}}
\newcommand{\Qabs}{Q_{\operatorname{abs}}}
\newcommand{\Qdirectvec}{\m{Q}_{\operatorname{direct}}}
\newcommand{\Qdirect}{Q_{\operatorname{direct}}}
\newcommand{\Qradvec}{\m{Q}_{\operatorname{rad}}}
\newcommand{\Qrad}{Q_{\operatorname{rad}}}
\newcommand{\Qreflvec}{\m{Q}_{\operatorname{refl}}}
\newcommand{\Qrefl}{Q_{\operatorname{refl}}}
\newcommand{\Qvis}{Q_{\operatorname{vis}}}
\newcommand{\Rsun}{R_{\odot}}
\newcommand{\Set}[1]{\left\{#1\right\}}
\newcommand{\Ssun}{S_{\odot}}
\newcommand{\Teq}{T_{\operatorname{eq}}}
\newcommand{\calI}{\mathcal{I}}
\newcommand{\calJ}{\mathcal{J}}
\newcommand{\esun}{e_{\odot}}

\newcommand{\xmax}{x_{\max}}
\newcommand{\xmid}{x_{\operatorname{mid}}}
\newcommand{\xmin}{x_{\min}}
\newcommand{\xsun}{\vec{x}_{\odot}}
\newcommand{\ymax}{y_{\max}}
\newcommand{\ymid}{y_{\operatorname{mid}}}
\newcommand{\ymin}{y_{\min}}

\graphicspath{{.}{Figures}}

\DeclareMathOperator{\diag}{diag}

\newcommand{\m}[1]{\boldsymbol{#1}}
\renewcommand{\vec}[1]{\m{#1}}

\newcounter{algocounter}
\renewcommand{\thealgocounter}{\arabic{algocounter}}
\newenvironment{algo}[1]{
  \refstepcounter{algocounter}
  \begin{mdframed}
    \noindent \textbf{Algorithm~\thealgocounter}: #1.
}{
  \end{mdframed}
}

\begin{document}

\begin{frontmatter}

\title{Fast hierarchical low-rank view factor matrices for thermal irradiance on planetary surfaces}

\author[1]{Samuel F. Potter\corref{cor1}}
\address[1]{Courant Institute of Mathematical Sciences, New York University, New York, NY, 10012}
\ead{sfp@cims.nyu.edu}

\author[2,3]{Stefano Bertone}
\address[2]{University of Maryland Baltimore County, Baltimore, MD 21250}
\address[3]{NASA Goddard Space Flight Center, Greenbelt, MD 20771}

\author[4]{Norbert Sch\"{o}rghofer}
\address[4]{Planetary Science Institute, Tucson, AZ 85719}

\author[3]{Erwan Mazarico}

\cortext[cor1]{Corresponding author}


\begin{abstract}
  We present an algorithm for compressing the radiosity view factor model commonly used in radiation heat transfer and computer graphics. We use a format inspired by  the hierarchical off-diagonal low rank format, where elements are recursively partitioned using a quadtree or octree and blocks are compressed using a sparse singular value decomposition---the hierarchical matrix is assembled using dynamic programming. The motivating application is time-dependent thermal modeling on vast planetary surfaces, with a focus on permanently shadowed craters which receive energy through indirect irradiance. In this setting, shape models are comprised of a large number of triangular facets which conform to a rough surface. At each time step, a quadratic number of triangle-to-triangle scattered fluxes must be summed; that is, as the sun moves through the sky, we must solve the same view factor system of equations for a potentially unlimited number of time-varying righthand sides. We first conduct numerical experiments with a synthetic spherical cap-shaped crater, where the equilibrium temperature is analytically available. We also test our implementation with triangle meshes of planetary surfaces derived from digital elevation models recovered by orbiting spacecrafts. Our results indicate that the compressed view factor matrix can be assembled in quadratic time, which is comparable to the time it takes to assemble the full view matrix itself. Memory requirements during assembly are reduced by a large factor. Finally, for a range of compression tolerances, the size of the compressed view factor matrix and the speed of the resulting matrix vector product both scale linearly (as opposed to quadratically for the full matrix), resulting in orders of magnitude savings in processing time and memory space.
\end{abstract}


\end{frontmatter}

\section{Introduction}

Illumination and temperature on the surfaces of airless bodies are dominated by radiative fluxes, from the Sun or from energy reflected or emitted from nearby topography. With the availability of high-resolution topographic maps of planetary surfaces, the modelling of radiative fluxes has become a major computational challenge. When a facet is irradiated by $N$ other facets within a direct line of sight, the total flux (radiance) received is the sum of $N$ fluxes. For example, Haworth Crater on the Moon ($86.9^\circ$S, $4^\circ$W) has a diameter of about 51~km and the laser altimeter-derived digital elevation model (DEM) currently available there has a horizontal resolution of 10~m/pixel. This means that fluxes from approximately $2\times 10^7$ facets must be summed to determine the incoming flux at each point for the entire crater surface, thus requiring on the order of $10^{14}$ FLOPs. This rough estimate includes only the first order of reflections and for a single position of the sun.

To dramatically reduce these computational costs, we present our implementation of a novel ``fast'' algorithm for this problem, with a focus on the rough surfaces of airless bodies such as planets and asteroids or features thereof. In computer graphics, global illumination algorithms which model light propagating in an environment with Lambertian surface reflectances go under the name ``radiosity''; this is essentially the model we consider here. Numerous fast radiosity solvers have been proposed in the computer graphics literature, but the calculation of the incident energy on planetary surfaces differs from these applications in several respects.

A major difference between planetary surfaces and the 3D models which were the focus of early radiosity algorithms is that planetary surfaces are discretized into a large number of small polygonal facets conforming to a rough terrain; on the other hand, early algorithms in computer graphics focused on contouring discontinuities in the solution of the radiosity integral equation under the assumption that the domain was composed of a relatively small number of large polygonal facets (e.g., walls, floors, ceilings). Contrary to the approaches taken in computer graphics, we compute a compressed version of the view factor matrix by assuming that hierarchically generated blocks of the radiosity kernel matrix can be compressed using low rank matrix factorizations. The standard form of a discretized radiosity equation is
\begin{equation}\label{eq:radiosity-system}
  \m{B} = \m{E} + \diag(\m{\rho}) \m{F} \m{B},
\end{equation}
where $\m{B}, \m{E}, \m{\rho} \in \mathbb{R}^N$ are vectors, $\m{F} \in \mathbb{R}^{N \times N}$ is a kernel matrix referred to as the ``view factor matrix'', and $N$ is the number of triangles. The rows and columns of $\m{F}$ correspond to individual triangular facets of a mesh conforming to a planetary surface, and $\m{F}_{ij}$ is an approximation of the view factor from the $j$th to $i$th triangle. Our assumption is that for two sets of triangles, indexed by $\calI, \calJ \subseteq \{1, \hdots, N\}$, the subblock $\m{F}_{\calI,\calJ}$ is ``numerically low rank''. That is, the singular values of $\m{F}_{\calI,\calJ}$ decay rapidly, and a low-rank approximation of $\m{F}_{\calI,\calJ}$ is available via the singular value decomposition. This is a reasonable assumption, since the radiosity system is the result of discretizing a second kind boundary integral equation. For comparison, with boundary integral methods for solving elliptic PDEs, the hierarchical matrix format (\Hmatrix)---different from what is proposed in this work---leads to an $O(N \log N)$ or better matrix vector product (MVP) with the kernel matrix in time and space~\cite{martinsson2019fast}. We demonstrate similar performance improvements for our format applied to $\m{F}$---see \Cref{sec:spherical-cap-crater,sec:gerlache}.

The radiosity equation is purely geometric. Provided that the surface does not change between time steps, the view factor matrix does not need to be recomputed. In the long-running, high-resolution, time-dependent thermal simulations needed in planetary science, a large number of radiosity MVPs with varying righthand sides are required. This is the justification for our approach: we pay an up front cost which is amortized over a large number of steps to gain access to fast MVPs at each time step.

The scientific applications of such an implementation are manifold. Remotely sensed temperatures are used to infer the thermal properties of the ground and albedo measurements within shadowed regions are affected by scattered light. The Moon has permanently shadowed regions (PSRs), which receive only indirect sunlight and their temperatures are determined by terrain irradiance. These extraordinarily cold regions are thought to harbor water ice~\cite{arnold79}, and are the targets of many upcoming landed and orbital missions to the Moon. The need to predict and quantify surface temperatures inside lunar PSRs is the main motivation behind the work presented here, but the same approach and implementation can be easily adapted to other planetary bodies and features.

Many of the currently available thermal models for planetary surfaces do not use fast algorithms for the evaluation of the terrain irradiance. Instead, they directly evaluate the scattered irradiance between each pair of facets, e.g.~\cite{mazarico11,delbo15,mazarico18,rubanenko2018ice,glaser19,schorghofer19,king20,hayne2021micro,mahanti22}. Heretofore, the most powerful model published to date for surface temperatures in lunar craters  is due to Paige et al.\ \cite{paige10a}. It achieves fast terrain irradiance calculations by constructing a mesh around a region of interest that is designed to optimize the radiance calculations. A fixed number of rays are cast from each pixel to determine visibility and radiances. The model is tuned to ensure agreement with measured temperatures, and is based on stochastic subsampling of the visible sky.

``Radiosity'' refers to a family of computer graphics algorithms for global illumination which solve the integral equation for Lambertian flux using the boundary element method~\cite{cohen1993radiosity}.  The term ``radiosity'' is overloaded: it refers to both this algorithm and also the radiative flux leaving a surface with units of W/m$^2$. This family of algorithms has mostly been displaced by ray tracing as the method of choice for global illumination. Nevertheless, the literature on this subject remains a useful source of techniques for solving other physical problems which involve the same integral equations, such as thermal radiation and particle flux. For example, radiosity was recently applied to free molecular flow~\cite{araki2020radiosity}.

The prototypical ``fast algorithm'' for radiosity, termed ``hierarchical radiosity'' is due to Hanrahan et al.~\cite{hanrahan91}.  Wavelet radiosity is a later generalization of hierarchical radiosity~\cite{gortler1993wavelet,kahler2008h}. These methods have all focused on accelerating radiosity in the particular case of computer graphics (at least as it was construed in the 1990s). In this case, to accurately contour the discontinuities in the solution (due to occlusion), hierarchical refinement and other methods were later used to accelerate the radiosity algorithm. Another line of research applied panel clustering, claiming $O(N \log^5 N)$ time in unoccluded environments and later $O(N^{3/2})$ time in occluded environments~\cite{atkinson2000numerical}. A separate line of research concerns Monte Carlo radiosity algorithms for computer graphics applications, usually referred to as ``stochastic radiosity''~\cite{bekaert2001hierarchical}. More recently, there have been some preliminary results related to using low-rank matrix factorizations for radiosity~\cite{fernandez09,aguerre2016hierarchical}.

\section{Physical description and governing equations}\label{sec:physical-model}

The surface energy balance on an airless body consists of direct solar irradiance, terrain irradiance, and subsurface heat flux.
The direct solar irradiance $\Qdirect$ is also known as insolation (incoming solar radiation). Its calculation requires knowledge of occlusion due to topography, which can block direct sunlight (terrain shadowing).
Terrain irradiance consists of a short-wavelength contribution $\Qrefl$ (sunlight reflected from land within view) and a long-wavelength contribution $\QIR$ (thermal infrared emitted from land within the field of view). The latter is also known as ``self-heating''.

On a surface $S$, the equation governing the energy balance is
\begin{equation}
(1 - \alpha) (\Qdirect + \Qrefl) + k \m{n}\cdot\nabla{T} + \epsilon \QIR = \epsilon \sigma T^4,
\label{eq:3d-surf}
\end{equation}
where $\alpha$ is the albedo, $\Qdirect$ is the direct incoming solar radiation, $k$ is the thermal conductivity, $T$ is the temperature, $\m{n}$ is the unit surface normal, $\epsilon$ is the emissivity, and $\sigma$ is the Stefan-Boltzmann constant.
In our work, we consider a ``Lambert albedo'', meaning that it is independent of the angle of incidence.
We will also assume that the subsurface heat flux is one-dimensional, as vertical heat fluxes are typically much larger than lateral heat fluxes on the large surfaces which we consider. Without subsurface heat conduction, the surface energy balance is
\begin{equation}
(1 - \alpha) (\Qdirect + \Qrefl) + \epsilon \QIR = \epsilon \sigma \Teq^4,
\label{eq:3d-equilbr}
\end{equation}
where $\Teq$ is in instantaneous equilibrium with the incoming flux and therefore known as the ``equilibrium temperature''.
The left-hand side is the absorbed flux.

On an airless body, if we denote the position of the Sun by $\xsun$, the insolation at a point $\vec{x}$ is
\begin{equation}
  \Qdirect = \frac{\Ssun}{\Rsun^2} \sin(\esun) V(\vec{x}, \xsun),\label{eq:insolation-for-point-sun}
\end{equation}
where $\Ssun$ is the solar constant, $\Rsun$ is the distance from the Sun in AU, and $\esun$ is the elevation of the sun, defined to be the angle between the tangent plane at $\vec{x}$ and the unit vector pointing from $\vec{x}$ to $\xsun$. The function $V$ is the line-of-sight visibility defined as
\begin{equation}\label{eq:visibility-function}
    V(\vec{x}, \vec{y}) = \begin{cases}
    1 & \mbox{if $\vec{y}$ is line-of-sight visible from $\m{x}$}, \\
    0 & \mbox{otherwise}.
    \end{cases}
\end{equation}
When the Sun is treated as a disk of finite extent rather than a point source, the flux can be expressed as a sum over point sources, the contribution of which takes the form of \eqref{eq:insolation-for-point-sun}.

The incoming radiances satisfy the two integral equations
\begin{align}
    \Qrefl(\vec{x}) &= \int_S F(\vec{x}, \vec{y}) \alpha(\vec{y}) \Big[\Qdirect(\vec{y}) + \Qrefl(\vec{y})\Big] dA(\vec{y}),
    \label{eq:Qrefl-integral-equation} \\
    \QIR(\vec{x}) &= \int_S F(\vec{x}, \vec{y}) \Big[\epsilon(\vec{y}) \sigma T(\vec{y})^4 + (1 - \epsilon(\vec{y})) \QIR(\vec{y})\Big] dA(\vec{y}), \label{eq:QIR-integral-equation}
\end{align}
where the view factor (sometimes also called the form factor) $F(\vec{x}, \vec{y})$ is given by
\begin{equation}
    F(\vec{x}, \vec{y}) = \frac{\cos\big(\theta(\vec{x}, \vec{y})\big) \cos\big(\theta(\vec{y}, \vec{x})\big)}{\pi |\vec{x} - \vec{y}|^2} V(\vec{x}, \vec{y}),
    \label{e:viewfactor}
\end{equation}
and where $\theta(\vec{x}, \vec{y})$ is the angle between the unit surface normal at $\vec{x}$ and the unit vector pointing from $\vec{x}$ to $\vec{y}$. The visibility function $V$ is defined by \eqref{eq:visibility-function}.

For comparison, the radiosity equation is the integral equation
\begin{equation}\label{eq:radiosity-integral-equation}
    B(\m{x}) = E(\m{x}) + \rho(\m{x}) \int_S F(\m{x}, \m{y}) B(\m{y}) dA(\m{y}),
\end{equation}
where $B$ is the radiosity (that is, the total radiant flux leaving $S$ per unit area), and $E$ is the directly emitted energy. This is the form used in computer graphics. The assumption of Lambertian reflectance implies~\cite{cohen1993radiosity}
\begin{equation}
  \alpha(\m{x}) Q(\m{x}) = B(\m{x}), \qquad \m{x} \in S ,
\end{equation}
so that the integral equations given by \eqref{eq:Qrefl-integral-equation} and \eqref{eq:QIR-integral-equation} are simply related to \eqref{eq:radiosity-integral-equation}. In terms of the reflected short-wavelength flux $\Qvis = \Qrefl + \Qdirect$, and using \eqref{eq:3d-surf} to substitute for temperature, we get
\begin{align}
  \Qvis(\m{x}) &= \Qdirect(\m{x}) + \int_S F(\m{x}, \m{y}) \alpha(\m{y}) \Qvis(\m{y}) dA(\m{y}), \label{eq:Qvis} \\
    \QIR(\m{x}) &= \int_S F(\m{x}, \m{y}) \Big[\big(1 - \alpha(\m{y})\big) \Qvis(\m{y}) + H(\m{y})\Big] dA(\m{y}) + \int_S F(\m{x}, \m{y})\QIR(\m{y}) dA(\m{y}), \label{eq:QIR}
\end{align}
where $H = k\, \m{n}\cdot \nabla T$ is the heat flux from the subsurface.

To model the heat conduction, we use the heat equation
\begin{equation}\label{eq:heat-equation}
    \rho c \frac{\partial T}{\partial t} = \nabla \cdot (k \nabla T),
\end{equation}
where $T$ is the temperature, $\rho c$ is the volumetric heat capacity and $k$ is the thermal conductivity. If $\m{n}$ is the outward facing unit surface normal for $S$, then we can combine \eqref{eq:heat-equation} with the radiation boundary condition
\begin{equation}\label{eq:radiation-boundary-condition}
  \Qabs(\m{x}) + k \m{n}(\m{x}) \cdot \nabla T(\m{x}) = \epsilon \sigma T(\m{x})^4, \qquad \m{x} \in S,
\end{equation}
where the absorbed flux is $\Qabs = (1 - \alpha) (\Qdirect + \Qrefl) + \epsilon \QIR$. In a simulation, the  runtime will be dominated by solving the system described in this section to be spent on the integral equations governing the flux. Our focus is thus to accelerate the solution of integral equations of the form \eqref{eq:radiosity-integral-equation}.

\section{Numerical discretization and occlusion}

We assume that the surface $S$ is available as a triangle mesh with an orientable boundary. We let $N$ denote the number of facets in the mesh. We denote the area, centroid, and surface normal of the $i$th triangle by $A_i$, $\m{x}_i$, and $\m{n}_i$, respectively. We approximate each surface integral by approximating functions defined on $S$ with piecewise constant basis functions, integrating over the triangle mesh, and using collocation to solve the integral equations. For instance, to discretize \eqref{eq:Qrefl-integral-equation}, we  approximate the integral over the triangle using the midpoint rule to get
\begin{equation}
    \Qrefl(\m{x}_i) = \sum_{j=1}^N F(\m{x}_i, \m{y}_j) \alpha(\m{y}_j) \Big[\Qdirect(\m{y}_j) + \Qrefl(\m{y}_j)\Big] A_j, \qquad 1 \leq i \leq N.
\end{equation}
This results in an $N \times N$ system with matrix $\m{F}$ whose entries are given by
\begin{equation}\label{eq:FF-def}
    \m{F}_{i, j} = F(\m{x}_i, \m{x}_j) A_j = \frac{\Big[\m{n}_i \cdot (\m{x}_j - \m{x}_i)\Big] \Big[\m{n}_j \cdot (\m{x}_i - \m{x}_j) \Big]}{\pi |\m{x}_i - \m{x}_j|^4} V(\m{x}_i, \m{x}_j) A_j.
\end{equation}
Altogether, we write the discretized version of \eqref{eq:Qrefl-integral-equation} as
\begin{align}\label{eq:Qrefl-system}
    \Qreflvec &= \m{F} \diag(\m{\alpha}) \Big[\Qdirectvec + \Qreflvec\Big],
\end{align}
and similarly for \eqref{eq:QIR-integral-equation}. Here, $\operatorname{diag}(\m{\alpha})$ is the $N \times N$ diagonal matrix whose $i$th diagonal entry is $\alpha(\m{x}_i)$, and $\Qreflvec$ (resp., $\Qdirectvec$) is a vector in $\mathbb{R}^N$ whose $i$th entry is $\Qrefl(\m{x}_i)$ (resp., $\Qdirect(\m{x}_i)$). Sophisticated methods for approximating view factors can be used instead of just the midpoint rule to compute the entries of $\m{F}$ more accurately. We use the midpoint rule for conceptual simplicity. The methods presented in this work are compatible with more accurate view factor approximations.

\Cref{eq:Qrefl-system} can easily be solved by using standard relaxation methods such as the Jacobi iteration. These methods will typically converge to machine precision in $O(1)$ iterations---there is no need for preconditioning or more elaborate schemes, so this is what we do. The number of nonzero entries scales like $O(N^2)$ as the mesh is refined, although the visibility factor $V(\m{x}, \m{y})$ results in a large number of zero entries. This is easy to see: consider two sets of triangles which are pairwise mutually visible, and assume that they account for $pN$ of the total triangles, where $0 < p < 1$. If this triangle mesh is refined uniformly, and in such a manner that this set of triangles remains mutually visible, then there will always be at least $p^2 N^2$ entries of $\m{F}$ which are nonzero. Hence, even in a scene with a high degree of occlusion (such as on the surface of an asteroid), an iterative radiosity solver based on a direct discretization cannot be faster than $O(N^2)$, thus motivating the need for an asymptotically faster matrix multiply.

\subsection{Hierarchically compressing the view factor matrix}\label{ssec:H-matrix-format}

The matrices stemming from discretized boundary integral equations are frequently  ``hierarchically off-diagonal low-rank'' (HODLR). The HODLR format and other similar formats such as the hierarchical semiseparable (HSS) format~\cite{xia2010fast} and the $\mathcal{H}$-matrix format~\cite{hackbusch2000sparse} have found widespread use in compressing the kernel matrices stemming from discretized integral equations and are useful in the design of ``fast direct solvers''. The HODLR format is not applicable to radiosity kernel matrices because of its directionality. Instead, we use a rather unstructured hierarchical matrix format assembled using dynamic programming, which we generically refer to as an \Hmatrix. We explain the dynamic programming algorithm used to assemble the compressed version of $\m{F}$ in \Cref{sec:hierarchical-low-rank-radiosity}.

We briefly describe our hierarchical matrix format. Let $\Omega \subseteq \mathbb{R}^d$ ($d = 2, 3$)
be a domain, $\partial \Omega$ its boundary, and
$\m{K} \in \mathbb{R}^{N \times N}$ a matrix of the form
\begin{equation}
  \m{K}_{ij} = K(\vec{x}_i, \vec{x}_j), \qquad \vec{x}_i, \vec{x}_j \in \partial \Omega, \qquad 1 \leq i \leq N, \qquad 1 \leq j \leq N .
\end{equation}
Let $\calI_1 \cup \cdots \cup \calI_P = [N]$ and let $\calJ_1 \cup \cdots \cup \calJ_Q = [N]$ define row and column index sets which partition $[N]$. Then, there are row and column permutation matrices $\m{P}$ and $\m{Q}$ such that
\begin{equation}
  \m{K} = \m{P} \begin{bmatrix}
    \m{K}_{\calI_1,\calJ_1} & \cdots & \m{K}_{\calI_1,\calJ_Q} \\
    \vdots & \ddots & \vdots \\
    \m{K}_{\calI_P,\calJ_1} & \cdots & \m{K}_{\calI_P,\calJ_Q}
  \end{bmatrix} \m{Q}^\top.
\end{equation}\label{eq:H-matrix-one-level}
These index vectors and permutations can be generated straightforwardly using spatial data structures such as quadtrees and octrees (see \Cref{ssec:quadtree}). We make fairly weak assumptions on the permuted subblocks of $\m{K}$ in \Cref{eq:H-matrix-one-level}. By recursively structuring $\m{K}$ using a spatial data structure, we create opportunities for compressing subblocks using a combination of their sparsity and low rank compression. The off-diagonal blocks are \emph{not} necessarily numerically low rank (which is instead the assumption for the HODLR format), we are still able to empirically obtain $O(N)$ complexity in cases of pratical interest (see \Cref{sec:spherical-cap-crater,sec:gerlache}).

\subsection{Visibility and raytracing}

A key component of assembling $\m{F}$ is evaluating the visibility function $V$ for a pair of points $\vec{x}_i$ and $\vec{x}_j$ on the surface $S$. For a pair of points $\vec{x}_i$ and $\vec{x}_j$ on the surface $S$, we require a means of testing whether the open interval $(\vec{x}_i, \vec{x}_j) = \{(1 - t)\vec{x}_i + t\vec{x}_j : 0 < t < 1\}$ intersects the volume at any point. If $\vec{x}_i$ and $\vec{x}_j$ ($i \neq j$) are two triangle centroids, we do this by shooting a ray from $\vec{x}_i$ to $\vec{x}_j$ and checking whether the ray intersects the $j$th triangle. The number of nonzero entries of the matrix $\m{F}$ equals the number of pairs of triangles for which $V(\vec{x}_i, \vec{x}_j) = 1$. So, if we let $\m{d}_{ij} = (\m{x}_j - \m{x}_i)/|\m{x}_i - \m{x}_j|$, we implement $V(x, y)$ as:
\begin{equation}
    V_{ij} = \begin{cases}
        1 & \mbox{if $\m{n}_i \cdot \m{d}_{ij} > 0$ and $\m{n}_j \cdot \m{d}_{ji} > 0$ and the ray shot from $\m{x}_i$ hits triangle $j$} \\
        0 & \mbox{otherwise}
    \end{cases}
\end{equation}
If $\m{n}_i \cdot \m{d}_{ij} \leq 0$ or $\m{n}_j \cdot \m{d}_{ji} \leq 0$, then there is no need to shoot a ray to test for occlusion, in which case the cost of evaluating $V(\m{x}_i, \m{x}_j)$ is $O(1)$---this is a fast and simple form of visibility culling, which can significantly reduce the number of rays that need to be traced for the types of meshes considered.

We experimented with two libraries which provide acceleration structures for visibility testing: CGAL~\cite{fabri2009cgal} and Embree~\cite{wald2014embree}. CGAL is a library of highly accurate computational geometry algorithms. It provides an axis-aligned bounding box (AABB) tree which allows very accurate raytracing in double precision. Embree is a more recent library which focuses on high-performance raytracing for computer graphics and rendering. It currently only provides a single-precision bounding volume hierarchy (BVH). In practice we find CGAL is about 3--5 times slower. On the other hand, Embree only supports single-precision floating point numbers, while CGAL can handle arbitrary precision arithmetic if necessary. We treat these libraries as black boxes and assume that if we build an acceleration data structure on a triangle mesh with $N$ triangles, we can evaluate ray-triangle intersections in $O(\log N)$ time. For our numerical tests, we use CGAL.

\subsection{Computing the full view factor matrix}\label{ssec:full-view-factor-matrix}

In order to compute subblocks of the view factor matrix $\m{F}$, we use the following algorithm to assemble it in the compressed sparse row (CSR) format:
\begin{algo}{Assemble the full view factor matrix in CSR format}
  \label{algo:FF-csr}
  \begin{enumerate}
  \item Fix a subset of row indices $\calI \subseteq [N]$ and column indices $\calJ \subseteq [N]$.
  \item Iterate over each triangle $\m{x}_i$ ($i \in \calI$), and find the indices $j \in \calJ \backslash \{i\}$ for which $\m{n}_i \cdot \m{d}_{ij} > 0$ and $\m{n}_j \cdot \m{d}_{ji} > 0$. For each such index $j$, shoot a ray from $\m{x}_i$ to $\m{x}_j$ and test for occlusion.
  \item If the $i$th and $j$th triangles are mutually visible, emit the row index pointer, the column index, and the value of $F_{ij}$ to incrementally construct the CSR matrix.
  \end{enumerate}
\end{algo}
This algorithm reduces the number of ray-triangle intersections which need to be performed and only requires us to compute entries of $F_{ij}$ which are nonzero.

\section{Hierarchical low-rank radiosity}\label{sec:hierarchical-low-rank-radiosity}

Several fast hierarchical summation algorithms for radiosity have been developed (see, e.g., \cite{hanrahan91,gortler1993wavelet}). Our application differs in that we would like a fast $\m{F}$ multiply in the case where the underlying mesh consists of a large number of small, triangular facets which contour a rough surface. To this end, our method is inspired by more recent kernel-free fast algorithms based on numerical linear algebra. We assume that there exists some recursive permutation of the entries of $\m{F}$ which creates opportunities for low rank compression at each level of the recursion. Because of the directionality and sparsity of the kernel, we must continue to refine off-diagonal blocks, unlike more familiar formats (see \Cref{ssec:H-matrix-format}). Despite having a large, constant percentage of zeros, the number of nonzero entries of $\m{F}$ scales like $O(N^2)$. After $\m{F}$ is approximated using this hierarchical format, we can finally think of it as being ``data sparse''~\cite{martinsson2019fast}.

\subsection{Spatial data structures: quadtrees and octrees}\label{ssec:quadtree}

Finding the desired permutation of the entries of $\m{F}$ is straightforward. They can be generated using a spatial data structure such as a quadtree or octree. To explain how this is done, we consider the scattered flux inside an ideal large crater, whose depth is given as the graph of a function: $z = z(x, y)$, so that each triangle centroid has the form $p_i = (x_i, y_i, z(x_i, y_i))$. A quadtree is a tree whose internal nodes correspond to a subset of $\mathbb{R}^2$ and a partition thereof into at most four point sets, and whose leaf nodes correspond to point sets in $\mathbb{R}^2$. Each internal node contains all of the points in each of its descendant nodes. Quadtrees differ in how they are constructed, for instance in the partitioning scheme used and in how the recursion is terminated. We use a quadtree constructed as follows:
\begin{algo}{Build quadtree on triangle mesh}
  \label{algo:quadtree}
  \begin{enumerate}
  \item Initially set the roots of the tree to be the centroids of the triangle mesh.
  \item Let $\xmin = \min_{1 \leq i \leq N} \m{x}_i$ be the minimum $x$-coordinate of each centroid, and let $\xmax, \ymin$, and $\ymax$ be defined similarly.
  \item Let $\xmid = (\m{x}_{\min} + \m{x}_{\max})/2$ and likewise define $\ymid$.
  \item Partition the triangle centroids associated with the current node into four sets:
    \begin{align*}
      \calI_1 &= \Set{i : \m{x}_i < \xmin \mbox{ and } y_i < \ymin}, \\
      \calI_2 &= \Set{i : \m{x}_i < \xmin \mbox{ and } y_i \geq \ymin}, \\
      \calI_3 &= \Set{i : \m{x}_i \geq \xmin \mbox{ and } y_i < \ymin}, \\
      \calI_4 &= \Set{i : \m{x}_i \geq \xmin \mbox{ and } y_i \geq \ymin}.
    \end{align*}
  \item For each index set $\calI_1, \calI_2, \calI_3$, and $\calI_4$, if the index set is nonempty, construct an internal node and recursively partition the centroids, starting from Step 2.
  \end{enumerate}
\end{algo}
If $\calI$ is the index set for a particular internal node, and $\calJ$ is the index set of one of its children, then $\calJ \subseteq \calI$. An octree is the analogous spatial data structure in 3D---the main difference is that each internal node has at most 8 child nodes instead of 4. Octrees are useful for topography with overhangs, or for small bodies, such as asteroids and comets.

\subsection{Low-rank off-diagonal blocks}\label{ssec:low-rank-off-diagonal-blocks}

Let $\calI$ and $\calJ$ be two disjoint sets of triangle indices, and let $\varepsilon > 0$ be a compression tolerance. We now describe a straightforward algorithm which attempts to compute a low-rank approximation of $\m{F}_{\calI,\calJ}$, the corresponding subblock of $\m{F}$:
\begin{algo}{Estimate numerical rank of $\m{F}_{\calI,\calJ}$}\label{algo:estimate-rank}
  \begin{enumerate}
  \item Use the algorithm of \Cref{ssec:full-view-factor-matrix} to compute $\m{F}_{\calI,\calJ}$ in CSR format.
  \item Let $k$ be an initial guess for the numerical rank of $\m{F}_{\calI,\calJ}$ (e.g., $k \gets 8$).
  \item Compute the $k$-term truncated sparse SVD $\m{F}_{\calI,\calJ} = \m{U}_k \m{\Sigma}_k \m{V}_k^\top$ from the CSR matrix.
  \item Denote the $i$th singular value of $\m{F}_{\calI,\calJ}$ by $\sigma_i$ ($1 \leq i \leq k$).
  \item Let $q$ be the smallest integer such that $\sigma_{q+1}/\sigma_1 < \varepsilon$, if it exists. If no such $q$ exists, continue to the next step. Otherwise, check whether:
    \begin{equation}\label{eq:nbytes-comparison}
      \mathtt{nbytes}(\m{U}_q) + \mathtt{nbytes}(\m{\Sigma}_q) + \mathtt{nbytes}(\m{V}_q) < \mathtt{nbytes}(\m{F}_{\calI,\calJ}).
    \end{equation}
    If so, then the truncated SVD that we have found is a low-rank approximation of $\m{F}_{\calI,\calJ}$; otherwise, we signal that we have failed to find one. \emph{Note}: here, we look at that size of the $q$-truncated singular vectors. It is also critically important to store $\m{U}_q$ and $\m{V}_q$ in the CSR format, since they will inherit a significant number of nonzeros from $\m{F}_{\calI,\calJ}$.
  \item Increase $k$ (e.g., by setting $k \gets 2k$), and jump back to Step 3.
  \end{enumerate}
\end{algo}
In general, we expect to be able to compute a $k$-term truncated SVD of an $m \times n$ matrix in roughly $O(kmn)$ time. This means that Algorithm~\ref{algo:estimate-rank} computes the $k$-term truncated SVD with the optimal number of terms in $O(k m n \log k)$ time. We compute SVDs using ARPACK~\cite{lehoucq1998arpack}.

In Algorithm~\ref{algo:estimate-rank}, the SVD is ``compatible'' with $\m{F}$ in the following sense. Let's say that $\m{F}$ is zero except for a subset of its entries indexed by two (possibly overlapping) index sets $\calI$ and $\calJ$. Then the maximum rank of $\m{F}$ is bounded by $\max(|\calI|, |\calJ|)$, where $|\cdot|$ is the cardinality of a set. (This can easily happen: consider a triangle mesh which is the discretization of a sphere with outward-facing surface normal; e.g., a perfectly spherical planet. In this case, $\calI = \calJ = \{\}$ and $\m{F} = \m{0}$.) The SVD will locate the index sets $\calI$ and $\calJ$ automatically and the left and right singular vectors will be correspondingly sparse. That is, there is no artificial fill-in hampering the likelihood of \eqref{eq:nbytes-comparison} to take effect.

Since the left and right singular matrices $\m{U}$ and $\m{V}$ are orthogonal matrices, they will nearly always have negative entries. Numerical errors then easily result in small negative fluxes when using our approximate view factor matrix. In practice we set any negative flux to zero.

We briefly comment that accelerating the assembly of $\m{F}_{\calI,\calJ}$ using something like adaptive cross approximation~\cite{zhao2005adaptive} or the interpolative decomposition~\cite{martinsson2019fast} may not be effective. Both methods assume that the spectrum of the matrix being compressed decays rapidly. This is not the case for the view factor matrix. This is not to say that it is impossible to compute a low rank approximation of $\m{F}_{\calI,\calJ}$ in subquadratic time, but our belief is that doing so will involve exploiting more detailed knowledge of the radiosiy kernel. We leave this for future work.

\subsection{Computing the hierarchically compressed view factor matrix}\label{ssec:assembly}

Using Algorithms~\ref{algo:FF-csr}, \ref{algo:quadtree}, and \ref{algo:estimate-rank}, we can assemble the full, hierarchically compressed approximation of $\m{F}$. With the compression tolerance $\varepsilon > 0$ fixed, we proceed as follows:
\begin{algo}{Assemble hierarchically compressed block $\m{F}_{\calI,\calJ}$}
  \label{algo:FF-assembly} \\
  \textbf{Input:} The index sets $\calI$ and $\calJ$ and a tolerance $\varepsilon > 0$. Optionally, the CSR version of the block $\m{F}_{\calI,\calJ}$. \\
  \textbf{Output:} A hierarchically compressed version of $\m{F}_{\calI,\calJ}$.
  \vspace{-0.25em}
  \begin{enumerate}
  \item Compute the block $\m{F}_{\calI,\calJ}$ using Algorithm~\ref{algo:FF-csr} if it was not supplied.
  \item If $|\calI| \cdot |\calJ|$ is less than a user defined threshold, set the current block to be either a sparse CSR or dense version of $\m{F}_{\calI,\calJ}$ depending on which uses fewer bytes.
  \item Otherwise, estimate the numerical rank $q$ of $\m{F}_{\calI,\calJ}$ using Algorithm~\ref{algo:estimate-rank} and compute the corresponding $q$-truncated sparse SVD.
  \item Recursively use this algorithm to hierarchically compress the subblocks of $\m{F}_{\calI,\calJ}$, taking care to pass along the corresponding subblock of the CSR version of $\m{F}_{\calI,\calJ}$.
  \item Compare the size of the sparse, dense, and truncated SVD version of $\m{F}_{\calI,\calJ}$ with the result of applying this algorithm recursively to each of its subblocks. Use the version of this block which requires the fewest number of bytes to be stored.
  \item Finally, if all of the blocks of this hierarchical matrix are sparse (or if all of them are dense), coalesce them into a single dense (or sparse) and set it to be the current block; otherwise, set the current block to be the hierarchical matrix.
  \end{enumerate}
\end{algo}
To compute $\m{F}$ itself, we apply Algorithm~\ref{algo:FF-assembly} to blocks at the first level of the spatial tree below the root (e.g., at most 16 pairs of index sets at the starting level of a quadtree, and 64 for an octree).

Getting an exact bound on the complexity of this algorithm is difficult, but in practice we observe Algorithm~\ref{algo:FF-assembly} to take a similar amount of time as it would take to compute the entire view factor matrix in the CSR format. When we run Algorithm~\ref{algo:FF-assembly}, we do not ever need to compute the entire CSR view factor matrix at once, only subblocks of it. For example, starting from the first level of a quadtree, where there are at most four blocks, we can work with the (at most) $4^2 = 16$ subblocks of the full view factor matrix one at a time. This signficantly reduces the amount of runtime memory use, which likely also has a positive impact on CPU runtime.

\subsection{Multiplying with the hierarchically compressed view factor matrix}

Let $\m{X}$ be a matrix with $N$ rows, and let $\m{X}_\calI$ be a matrix formed from a subset of its rows indexed by $\calI$. After assembling the hierarchically compressed view factor matrix using Algorithm~\ref{algo:FF-assembly}, we can easily compute the product $\m{Y} = \m{F}\m{X}$ using the following algorithm:
\begin{algo}{Compute matrix-matrix product with hierarchically compressed $\m{F}$}
  \label{algo:FF-multiply}
  \begin{enumerate}
  \item Initially set $\m{Y} \gets \m{0}$.
  \item Set the current node to the root of the quadtree (or octree).
  \item For each pair of index sets $\calI$ and $\calJ$ amongst $\calI_1, \calI_2, \calI_3$, and $\calI_4$:
    \begin{enumerate}
    \item If $\m{F}_{\calI, \calJ}$ is a CSR matrix, a dense matrix, or a truncated SVD, compute $\m{F}_{\calI, \calJ} \m{X}_\calJ$ in the obvious way.
    \item Otherwise, $\m{F}_{\calI, \calJ}$ is a block matrix. Recursively apply this algorithm to compute $\m{F}_{\calI, \calJ} \m{X}_\calJ$.
    \item Set $\m{Y}_\calI \gets \m{Y}_\calI + \m{F}_{\calI,\calJ} \m{X}_\calJ$.
    \end{enumerate}
  \end{enumerate}
\end{algo}
To compute a matrix-vector product, we can take $\m{X}$ to be a matrix with a single column.

HODLR-compressed matrices' MVPs have a time complexity of $O(N \log N)$. Because of the more complicated structure of the view factor matrix, it is not clear whether we can expect an asymptotically improved MVP in the general case, since it depends delicately on $\partial\Omega$. That said, we can expect it to do no worse than multiplying with the uncompressed view factor matrix because of the way we construct $\m{F}$ in Algorithm~\ref{algo:FF-assembly}. For a reasonable choice of parameters, we observe $O(N)$ scaling---see \Cref{sec:spherical-cap-crater,sec:gerlache}.

\section{Time integration}

\subsection{Equilibrium temperatures}\label{sec:time-independent}

To solve the time-independent problem, for a particular Sun position, we compute $\Qdirectvec$ by checking whether the Sun is visible from a centroid $\m{x}_i$ using raytracing and computing $\Qdirect(\m{x}_i)$ from \eqref{eq:insolation-for-point-sun}. We then rearrange Eq.~\eqref{eq:Qrefl-system} to get
\begin{equation}
    \Big(\m{I} - \m{F} \operatorname{diag}(\m{\alpha})\Big) \Qreflvec = \m{F} \operatorname{diag}(\m{\alpha}) \Qdirectvec. \label{eq:time-indep-system-1},
\end{equation}
which can be solved in a small number of iterations using a simple Jacobi iteration. Discretizing \eqref{eq:QIR} gives:
\begin{equation}
    \Big(\m{I} - \m{F}\Big) \QIRvec = \m{F} \big(\m{I} - \operatorname{diag}(\m{\alpha})\big) \big(\Qdirectvec + \Qreflvec\big) + \m{F}\m{H}. \label{eq:time-indep-system-2}
\end{equation}
This system can be solved as before, yielding $\QIRvec$. Finally, we can compute $\m{T}$ from Eq.~\eqref{eq:3d-equilbr}. Equations~\eqref{eq:time-indep-system-1} and~\eqref{eq:time-indep-system-2} are both of the form Eq.~\eqref{eq:radiosity-system}, with $\m{\rho} = \m{\alpha}$ for $\m{B} = \Qreflvec$ and $\m{\rho} = \m{1}$ for $\m{B} = \QIRvec$.

\subsection{Time integration with subsurface heat storage}

To solve the time-dependent problem, we only consider 1D temperature profiles beneath the centroid of each triangle. For each centroid $\m{x}_i$, we let $z_j > 0$ be the depth beneath the surface in the direction $\m{n}_i$ (i.e., the point $\m{x}_i - z_j \m{n}_i$), and let $M$ denote the number of layers, so that $1 \leq j \leq M$. We discretize \eqref{eq:heat-equation} in 1D with this set of grid points in a flux-conservative manner and use a Crank-Nicolson method with a linearized upper boundary condition. Details are described in~\cite{schorghofer2017planetary}.
We obtain a boundary condition for the surface from \eqref{eq:radiation-boundary-condition}. For a boundary condition in the interior, we prescribe a fixed heat flux $\m{H}$, which can be zero. We typically choose $M$ fairly small ($M$ is typically $20$ to $30$), independently of $N$ (the number of triangle faces). The cost of driving each subsurface thermal model at each step is then $O(NM)$ overall, since the time step only involves solving an $M \times M$ tridiagonal system.

To solve the time-dependent problem, we could take an approach similar to the one presented in \cref{sec:time-independent}, but we consider a more cost-effective alternative, reducing the cost of each time step while maintaining accuracy. Subsurface heat conductance introduces a history dependence in the problem. In this case, there is little incentive to solve the systems describing $\Qreflvec$ and $\QIRvec$ at each time step. Multiplying with $\m{F}$ once corresponds to computing one order of scattering. If we multiply once at each time step, many orders of reflections are taken into account, with the higher order reflections lagging behind in time. We use the following time-stepping scheme for the time-dependent problem:

\begin{algo}{Run time-dependent simulation}
  \label{algo:time-dependent}
  \\ \\
  Let $\Qabsvec^{(n)}$, $\Qreflvec^{(n)}$, $\QIRvec^{(n)}$, and $\m{T}^{(n)}$ be vectors giving the values of $\Qabs$, $\Qrefl$, $\QIR$, and $T$ over each triangle at the $n$th time step. To ease notation, we let $\Qradvec^{(n)}$ be defined by:
  \begin{equation}
    \Qrad(\m{x}_i) = \epsilon(\m{x}_i) \sigma T(\m{x}_i)^4, \qquad 1 \leq i \leq N.
  \end{equation}
  For $n = 0$, we initialize:
  \begin{align}
    \Qabsvec^{(0)} &= \big(\m{I} - \diag(\m{\alpha})\big) \Qdirectvec^{(0)} + \m{H}^{(0)}, \\
    \Qreflvec^{(0)} &= \m{0}, \\
    \QIRvec^{(0)} &= \m{0}, \\
    \Qradvec^{(0)} &= \Qabsvec^{(0)}.
  \end{align}
  Then, at each time step, we update the fluxes by:
  \begin{align}
    \Qreflvec^{(n+1)} &= \m{F} \operatorname{diag}(\m{\alpha}) \Big[\Qdirectvec^{(n+1)} + \Qreflvec^{(n)}\Big], \\
    \QIRvec^{(n+1)} &= \m{F} \Big[\Qradvec^{(n)} + \big(\m{I} - \diag(\m{\epsilon})\big) \QIRvec^{(n)}\Big], \\
    \Qabsvec^{(n+1)} &= \big(\m{I} - \diag(\m{\alpha})\big) (\Qdirectvec^{(n+1)} + \Qreflvec^{(n+1)}) + \diag(\m{\epsilon}) \QIRvec^{(n+1)}.
  \end{align}
  We then update the subsurface temperature profiles by stepping each Crank-Nicolson scheme, with $\Qabsvec^{(n)}$ and $\Qabsvec^{(n+1)}$ as input. This yields the surface temperature $\m{T}^{(n+1)}$, and hence $\Qradvec^{(n+1)}$.
\end{algo}

\section{Implementation}

In this work, we present a prototype implementation which proves the concept---in a future work, we will more carefully address its optimization.
The algorithms described here are available as a Python package:
\begin{center}
    \url{https://github.com/sampotter/python-flux}
\end{center}
Extensive use is made of numpy~\cite{van2011numpy} and scipy~\cite{virtanen2020scipy}. In particular, we nearly always store full view factor matrix blocks as instances of the \texttt{scipy.sparse.csr\_matrix} class. To compute truncated SVDs, we use \texttt{scipy.sparse.linalg.svds}, which is a wrapper around ARPACK~\cite{lehoucq1998arpack}. The combination of ARPACK and \texttt{csr\_matrix} allows us to compute truncated SVDs of view factor matrix blocks in time proportional to the number of nonzero entries, which is critically important for the efficiency of our implementation. Since a large amount of numpy is parallelized ``behind the scenes'', our implementation benefits from this, even though we do not focus on parallelism in this work.

\section{Numerical validation: equilibrium temperature of a spherical cap crater}\label{sec:spherical-cap-crater}

\begin{figure}
    \centering
    \begin{subfigure}[b]{0.45\textwidth}
        \includegraphics[width=\textwidth]{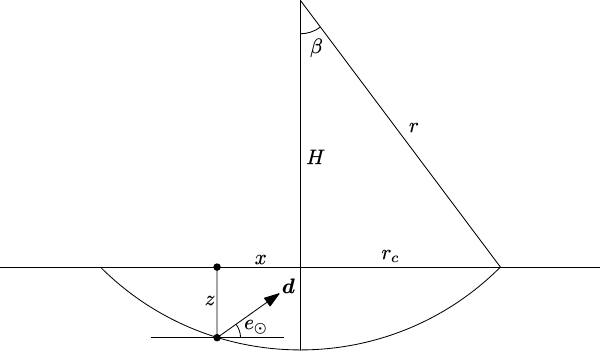}
        \caption{Cross section of a spherical cap-shaped crater, where
        $\beta$ is the angle between the crater rim and axis of the crater,
        $r$ is the radius of the sphere,
        $H$ is the distance from the sphere center to the ground plane,
        $r_c$ is the radius of the crater,
        $d_{\otimes}$ is the direction to the sun,
        and $\esun$ is elevation angle of the sun.
        }\label{fig:ingersoll-crater}
    \end{subfigure}
    \hfill
    \begin{subfigure}[b]{0.5\textwidth}
        \centering
        \includegraphics[width=\linewidth]{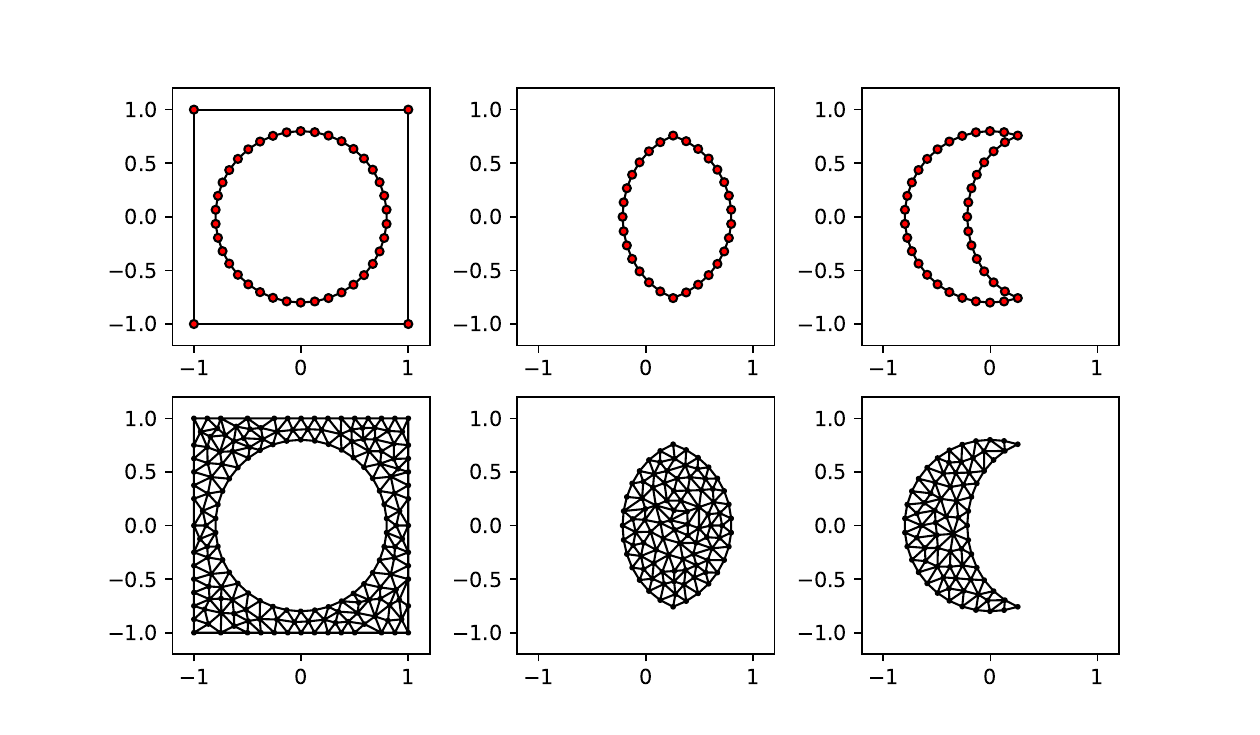}
        \caption{An example of meshing a hemispherical crater into three parts. The parameters for this problem are $h = {(2/3)}^5 = 0.131\ldots$, $e_0 = 15^\circ$, $\vec{d} = (\cos e_0, 0, \sin e_0)$, $\beta = 40^\circ$, and $r_c = 0.8$. \emph{Top}: the contours of each part of the mesh. \emph{Bottom}: the meshes, generated using Shewchuk's Triangle program~\cite{Shewchuk:1996aa}. \emph{Left}: the horizontal ground plane. \emph{Middle}: the shadowed part of the crater. \emph{Right}: the illuminated part of the crater.}
        \label{fig:ingersoll-meshing}
    \end{subfigure}
    \caption{Depiction of the spherical cap crater example.}
\end{figure}

\subsection{Analytical solution for a bowl-shaped crater}

We validate our implementation against the exact solution for the temperature, which is available in the ideal case of a spherical cap.
In this case, the equilibrium temperature in the shadow is given by \citep{buhl68,ingersoll92}
\begin{equation}
\sigma \Teq^4 = \Ssun \sin \esun \, \frac{f (1 - \alpha)}{1 - \alpha f}
\left[ 1+ \frac{\alpha (1 - f)}{\epsilon} \right] ,
\end{equation}
where $\alpha$ is the albedo, $\Ssun$ is the solar constant, $\esun$ the elevation of the Sun relative to the horizontal, and $\epsilon$ the infrared emissivity. These parameters match \Cref{sec:physical-model}. The geometric parameter $f$ is given by the diameter-to-depth ratio, $D/d$, of the crater:
\begin{equation}
\frac{1}{f} = 1+ \frac{1}{4} \left(\frac{D}{d}\right)^2
\end{equation}
The temperature in each part of the domain is given by the indirect irradiance, which is constant within the spherical portion, plus the direct irradiance, such that
\begin{equation}
\epsilon\sigma \Teq^4 =
(1 - \alpha) \Ssun \times \left\{ \begin{array}{ll}
\sin \esun & \mbox{sunlit, outside of crater} \\
\sin e + b \sin \esun & \mbox{sunlit, in crater} \\
b \sin \esun & \mbox{shadowed, in crater} \\
\end{array} \right.,
\label{eq:globalsolution}
\end{equation}
where $e$ is the elevation of the Sun relative to the local surface, and
\begin{equation}
    b = f \frac{\epsilon + \alpha(1 - f)}{1 - \alpha f}.
\end{equation}
See \Cref{fig:ingersoll-crater} for a diagram depicting these parameters.

\subsection{View factor computation}

In this work, our goal is to integrate triangle-to-triangle view factors. Although we use a simple midpoint rule for quadrature, we mention that there are a number of ways of doing so analytically---e.g., \cite{schroder1993closed,walton2002calculation,narayanaswamy2015analytic}. Combining analytic integration of triangle-to-triangle view factors in the near-field with our currect approach is likely to provide a very reasonable scheme. We emphasize that Algorithm~\ref{algo:FF-assembly} remains unchanged if other methods of computing view factors are substituted. Such a change would result in different (more accurate) view factor matrix entries, but that is all. We leave exploration of this issue for future work.

\subsection{Crater and shadow contouring}

The solution of the system described in \Cref{sec:physical-model} is piecewise continuous on $S$. For the spherical cap-shaped crater considered in this section, there will be a jump discontinuity in $\Teq$ across the boundaries separating the three regions of $S$. It is important to assess the effect of this discontinuity on the accuracy of the equilibrium temperature, $\Teq$. To this end, we consider two discretizations of the crater: one where the triangle mesh itself contours the rim of the crater and the shadow exactly, and one where only the crater rim is contoured.

For a given $y_s$ ($-r_c \leq y_s \leq r_c$), a silhouette point $(x_s, y_s)$ satisfies:
\begin{equation}
    \sqrt{r^2 - x_s^2 - y_s^2} + \tan(e_0) \left(\sqrt{r_c^2 - y_s^2} - x_s\right) = H.
\end{equation}
For a point strictly in the interior of the crater ($x_s^2 + y_s^2 < r_c^2$), there are two values of $x_s$ that satisfy this equation: one corresponding to the silhouette line, and another which is on the crater rim which we denote by $x_c$. As $x_s^2 + y_s^2 \to r_c^2$, these two solutions coalesce. To find $x_s$, we use a rootfinder. This breaks down near the crater rim; in this case, we extrapolate $x_s$ from nearby $(y_s, x_s(y_s))$ pairs. We denote the maximum absolute $y$-coordinate of a point on the silhouette by $y_p > 0$.

\subsection{Numerical results}

\begin{figure}[tbh!]
  \begin{tabular}{cc}
  (a) & (b)\\
    \includegraphics[width=0.5\linewidth]{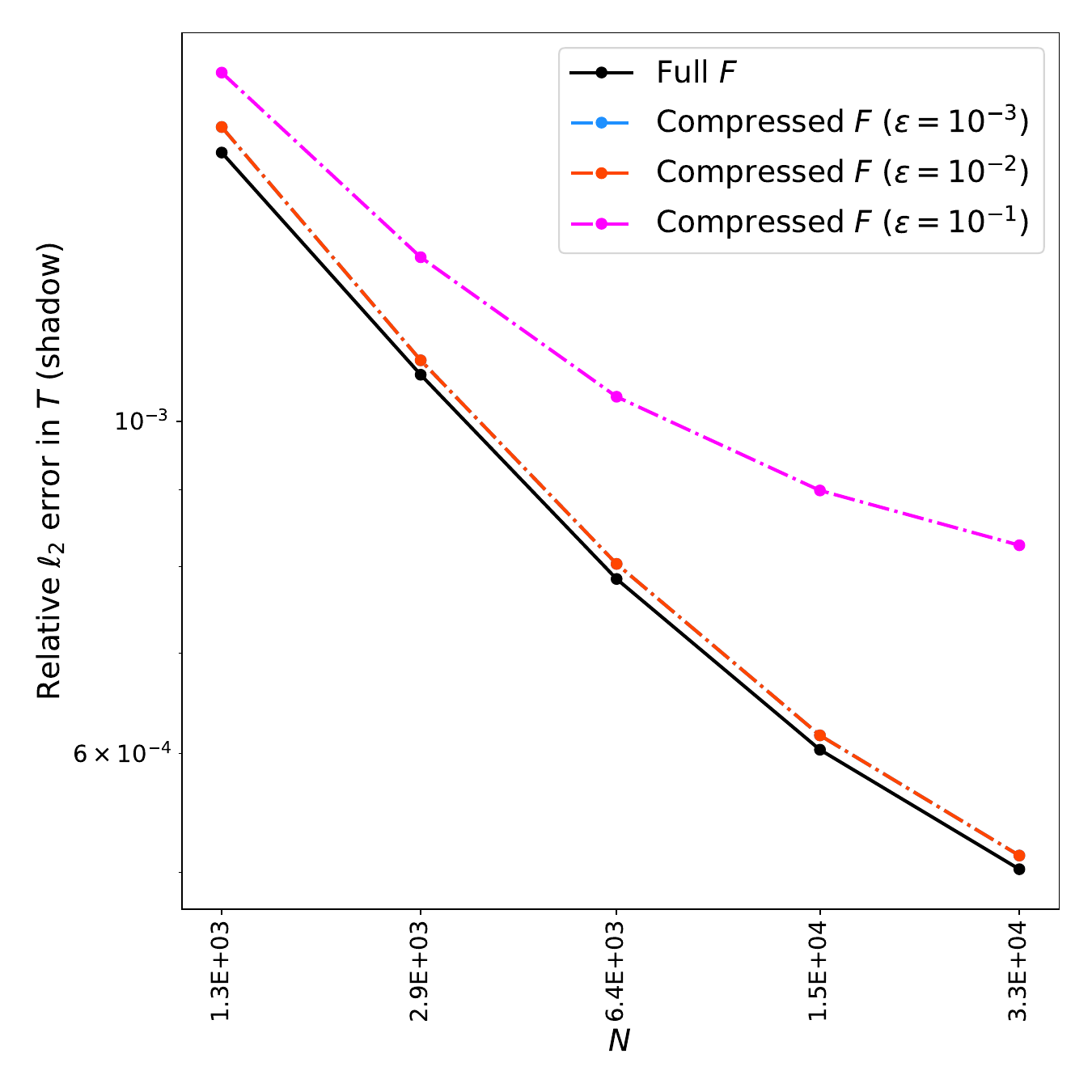}
  &
  \includegraphics[width=0.5\linewidth]{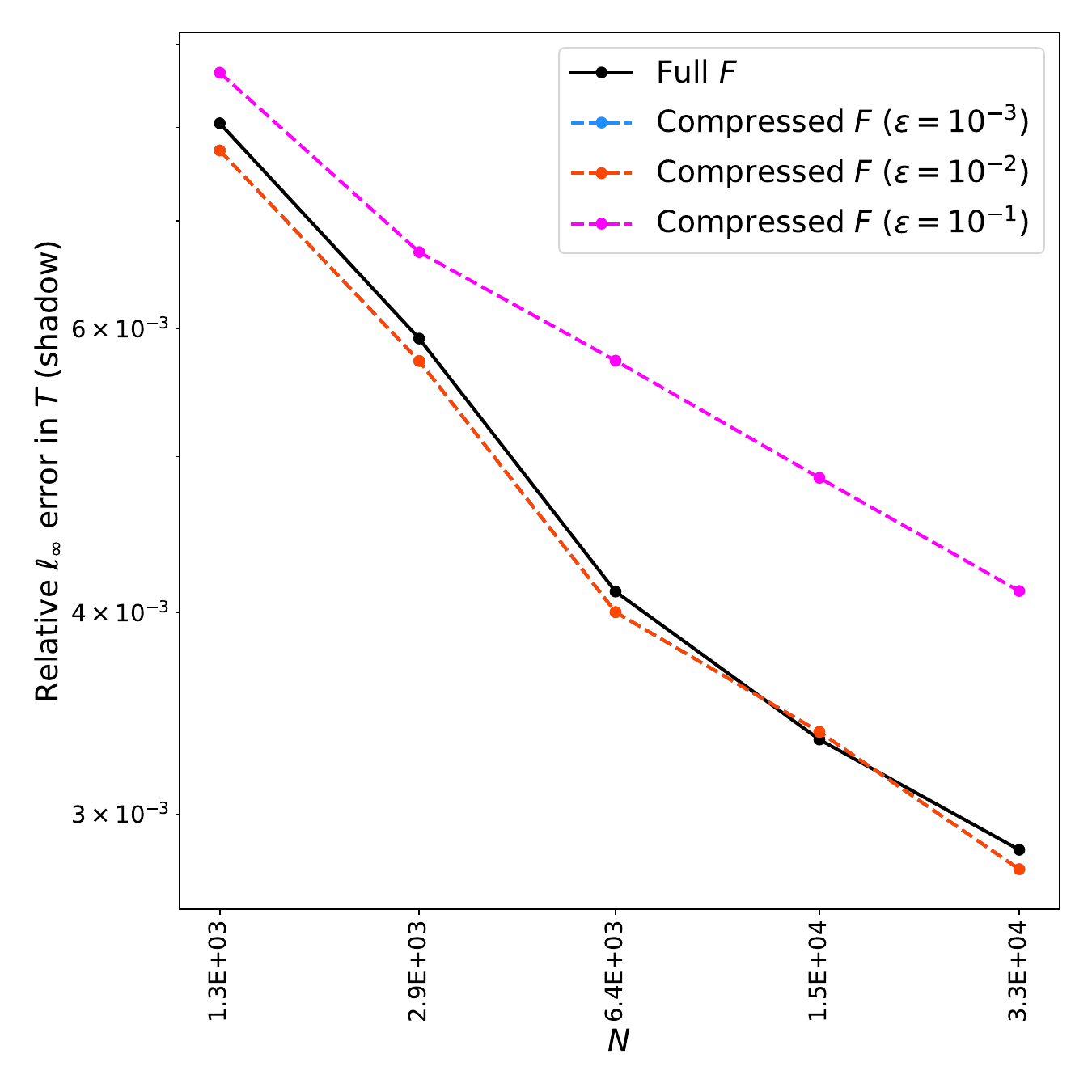} \\
  (c) & (d) \\
  \includegraphics[width=0.5\linewidth]{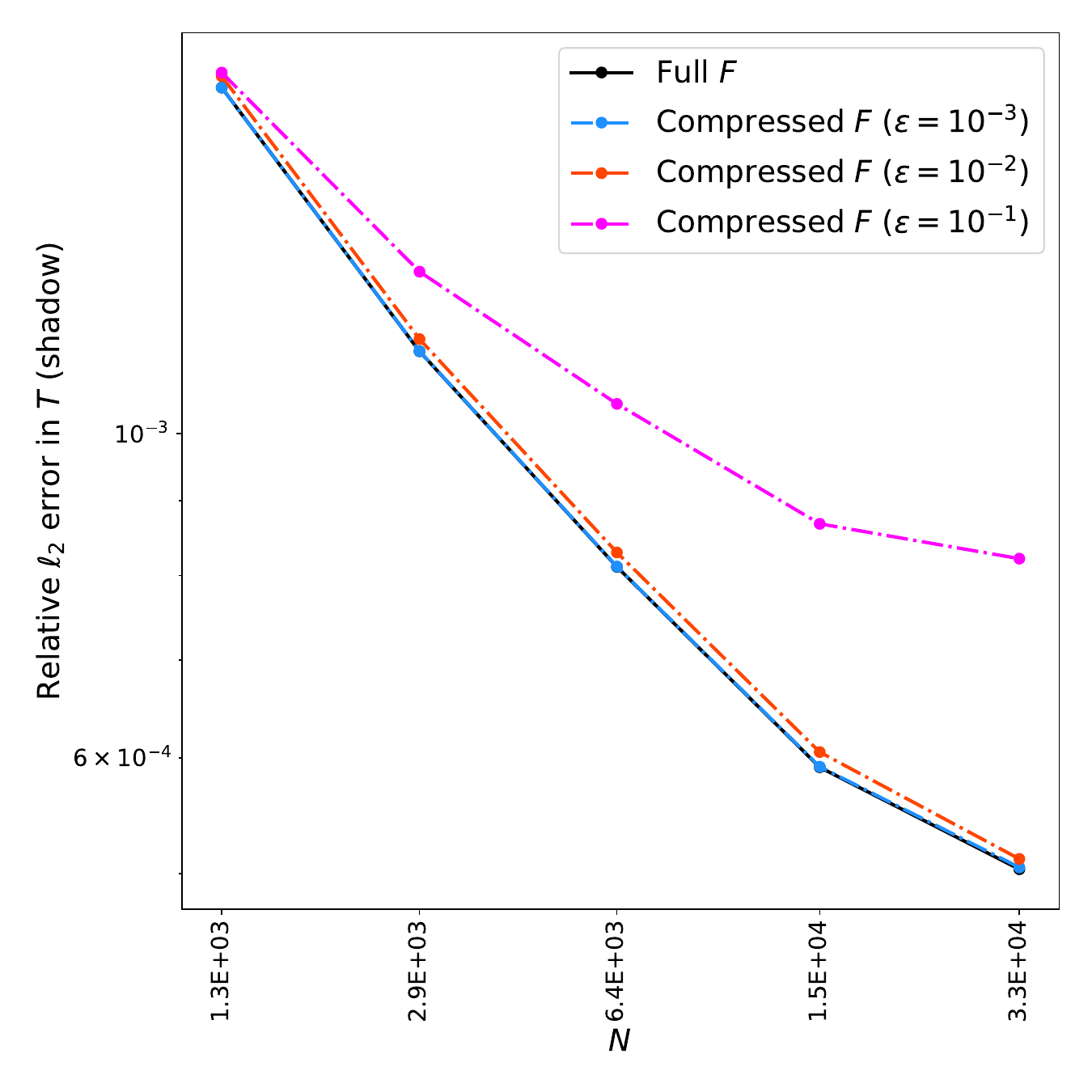}%
  &
  \includegraphics[width=0.5\linewidth]{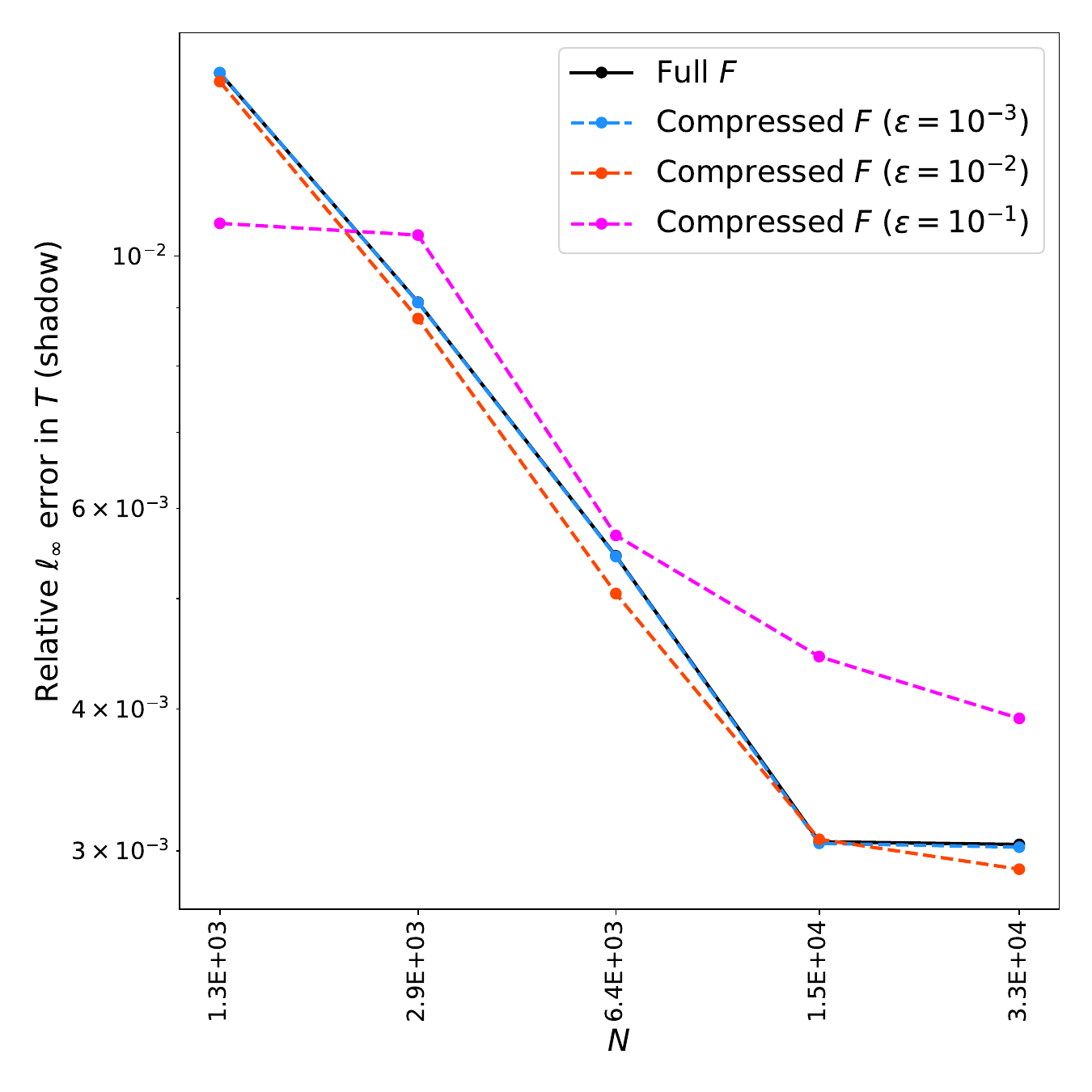}
  \end{tabular}
  \caption{Relative $\ell_1, \ell_2$, and $\ell_\infty$ errors for the spherical cap-shaped crater for varying $N$ (number of triangular facets). For $p = 1, 2, \infty$, $\|\hat{\m{T}} - \m{T}\|_p / \|\m{T}\|_p$ is computed, where $\m{T}$ is computed from \eqref{eq:globalsolution}, and $\hat{\m{T}}$ is computed using the discretized view factor matrix. Here, ``Full $\m{F}$'' refers to Algorithm~\ref{algo:FF-csr} and ``Compressed $\m{F}$'' refers to Algorithm~\ref{algo:FF-assembly}. \emph{Top row}: contouring both the crater rim and the shadow line. \emph{Bottom row}: contouring only the crater rim. \emph{Left column}: the relative $\ell_2$ error, where (a) uses full contouring while in (c) only crater rim is contoured. \emph{Right column}: the relative $\ell_\infty$ error, where (b) uses full contouring, and (d) only contours the crater rim.}\label{fig:spherical-crater-errors}
\end{figure}

\begin{figure}[tbh!]
  \begin{tabular}{cc}
  (a) & (b)\\
    \includegraphics[width=3.2in]{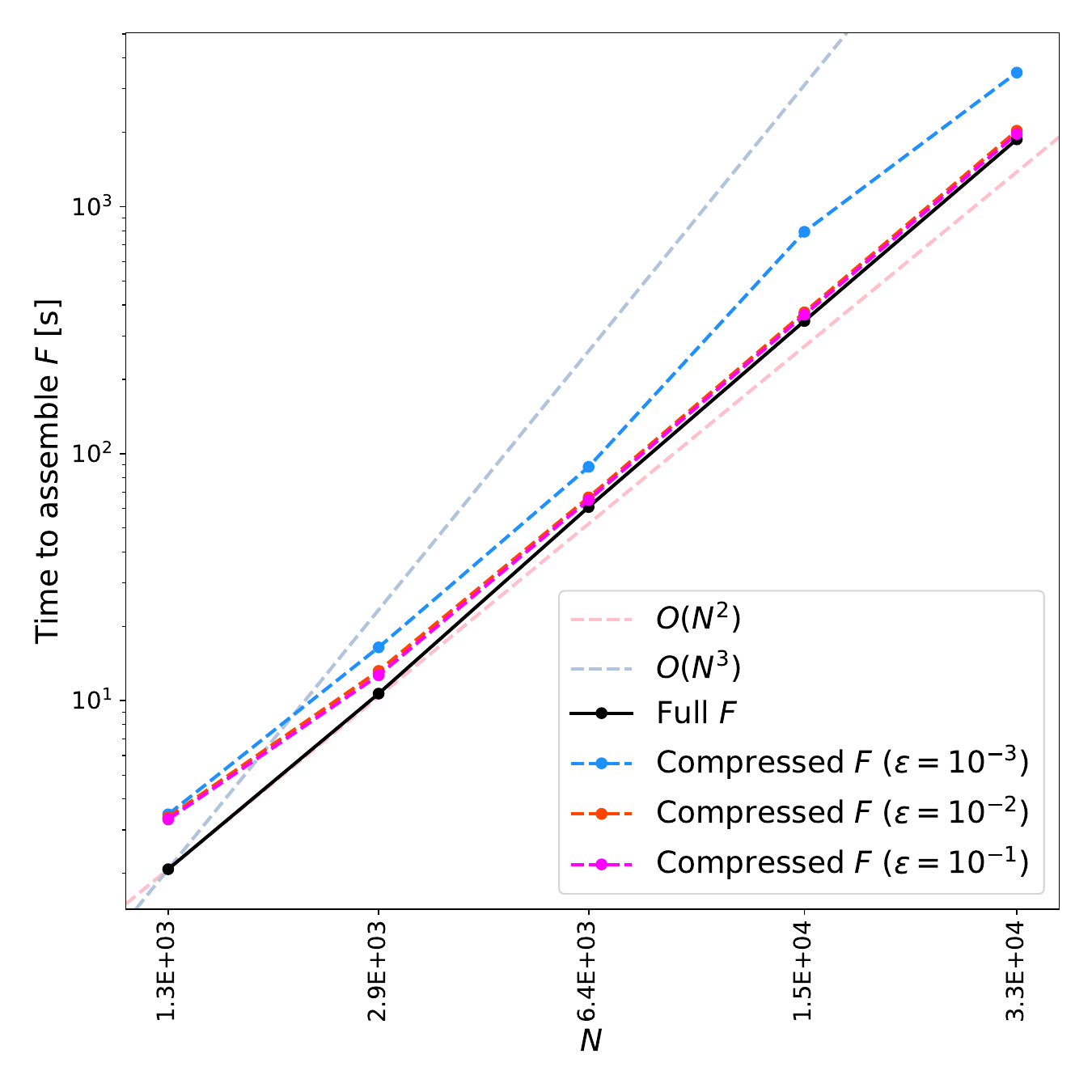}
  &
    \includegraphics[width=0.5\linewidth]{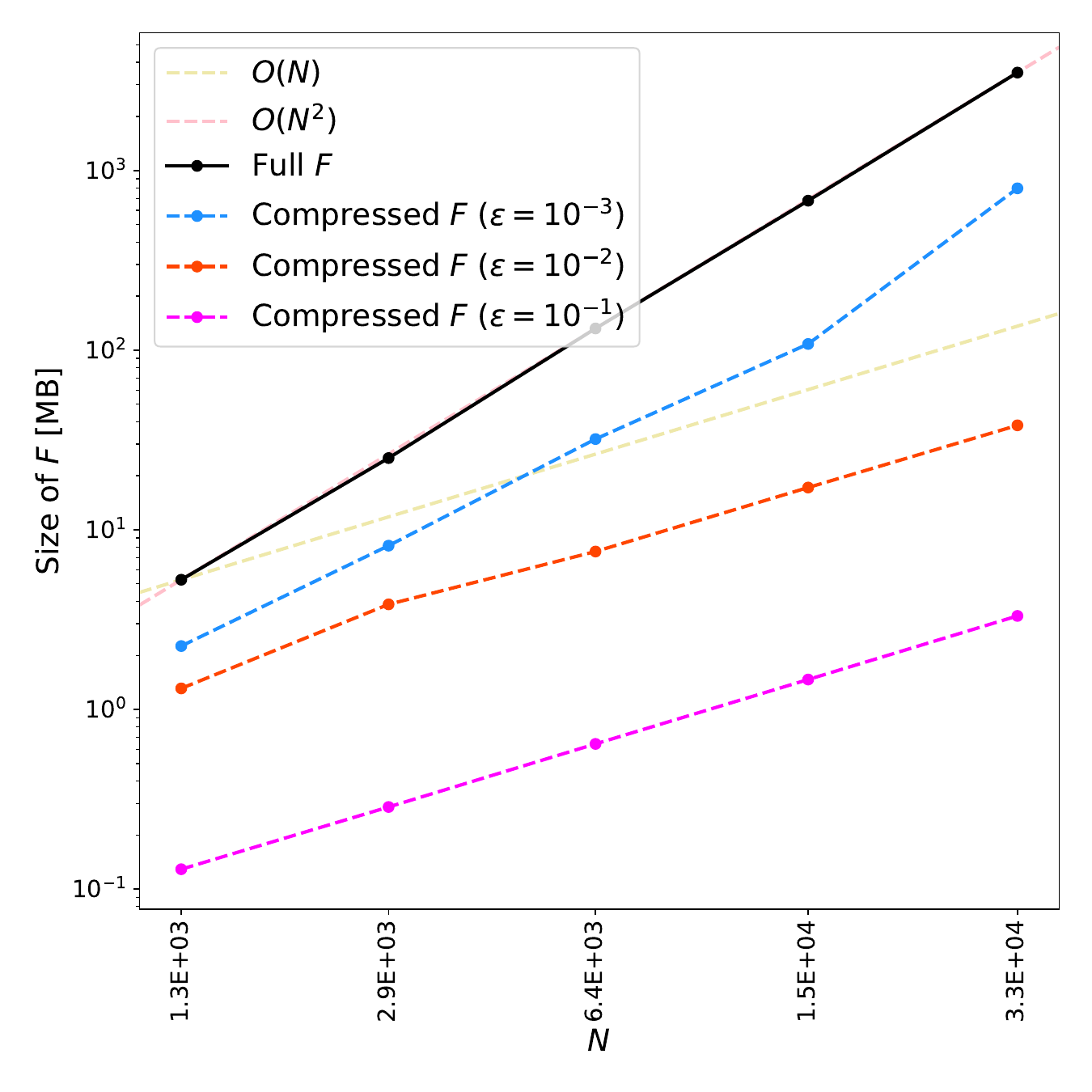} \\
  (c) & (d) \\
    \includegraphics[width=0.5\linewidth]{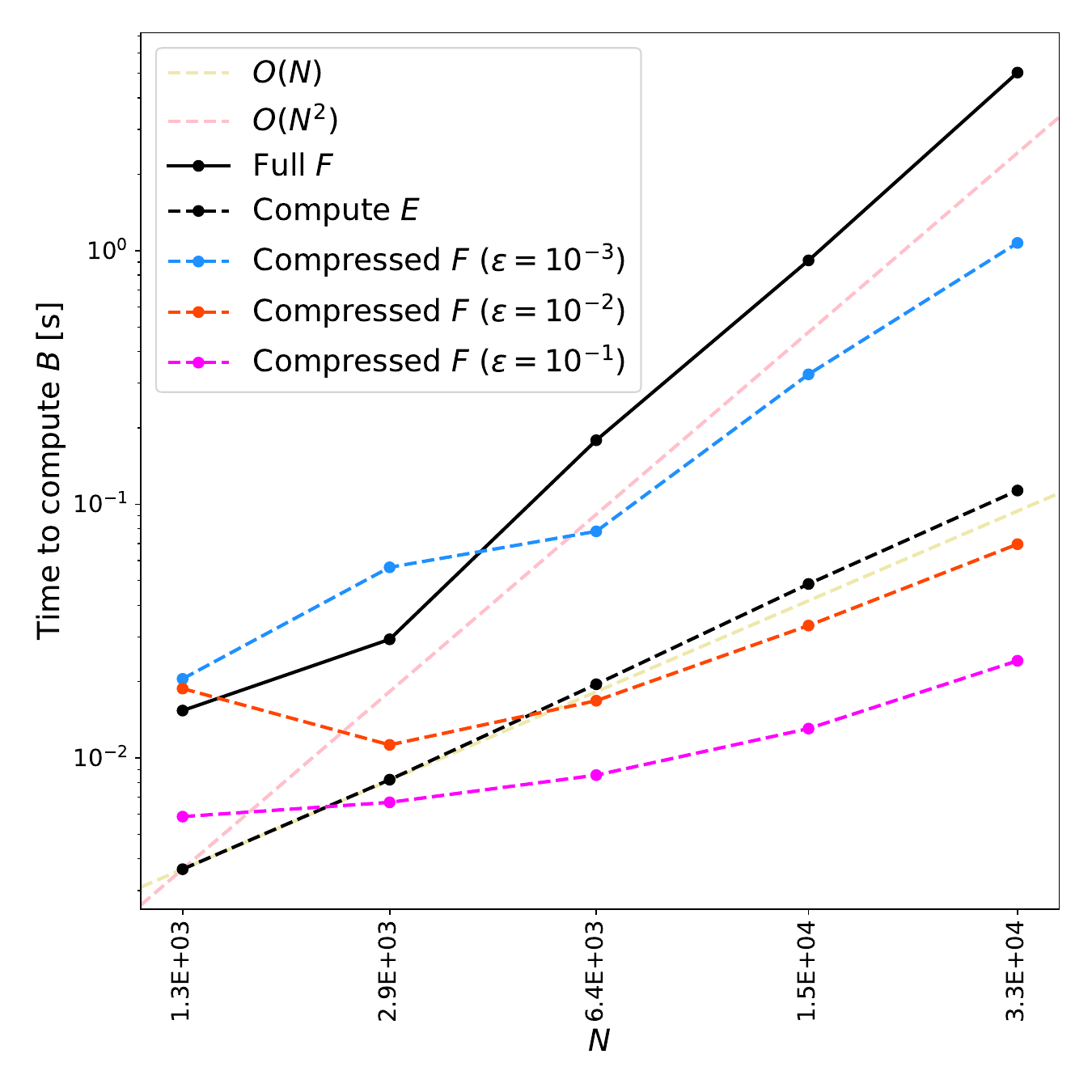}%
  &
  \includegraphics[width=0.5\linewidth]{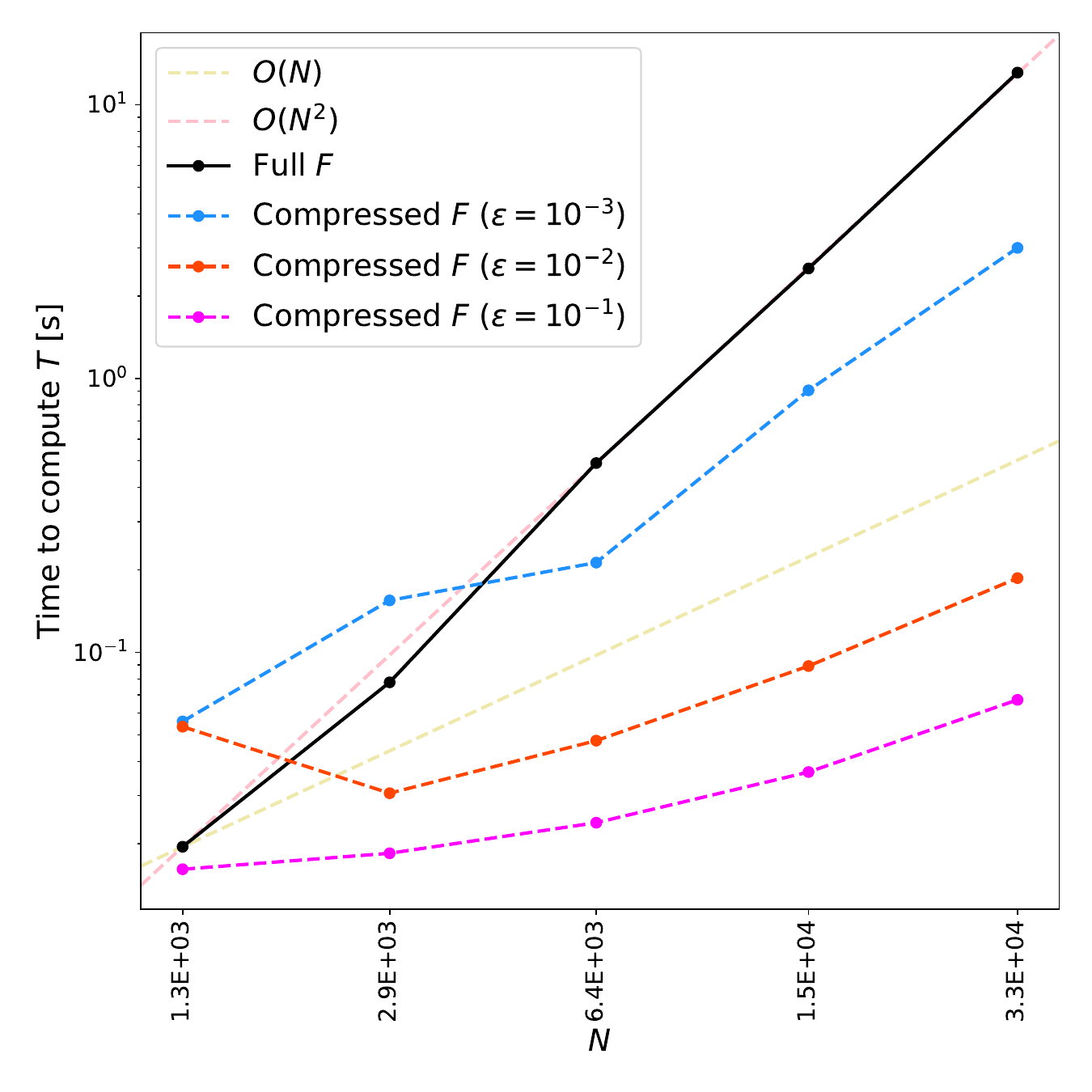}
  \end{tabular}
  \caption{Performance characteristics of the hierarchically compressed view factor matrix for the spherical cap-shaped crater as $N$ (the number of triangular facets) is varied. In each legend, ``Full $\m{F}$'' indicates Algorithm~\ref{algo:FF-csr} and ``Compressed $\m{F}$'' indicates Algorithm~\ref{algo:FF-assembly} for varying $\varepsilon$. In (a), we evaluate the time it takes ot assemble the full and compressed view factor matrices using Algorithm~\ref{algo:FF-assembly}; (b) shows the size in MB of the full and compressed view factor matrices; (c) shows the time it takes to compute $\m{B}$; and (d) shows the time it takes to compute $\m{T}$.}\label{fig:spherical-crater-performance-characteristics}
\end{figure}

\begin{figure}
  \centering
  \includegraphics[width=0.7\linewidth]{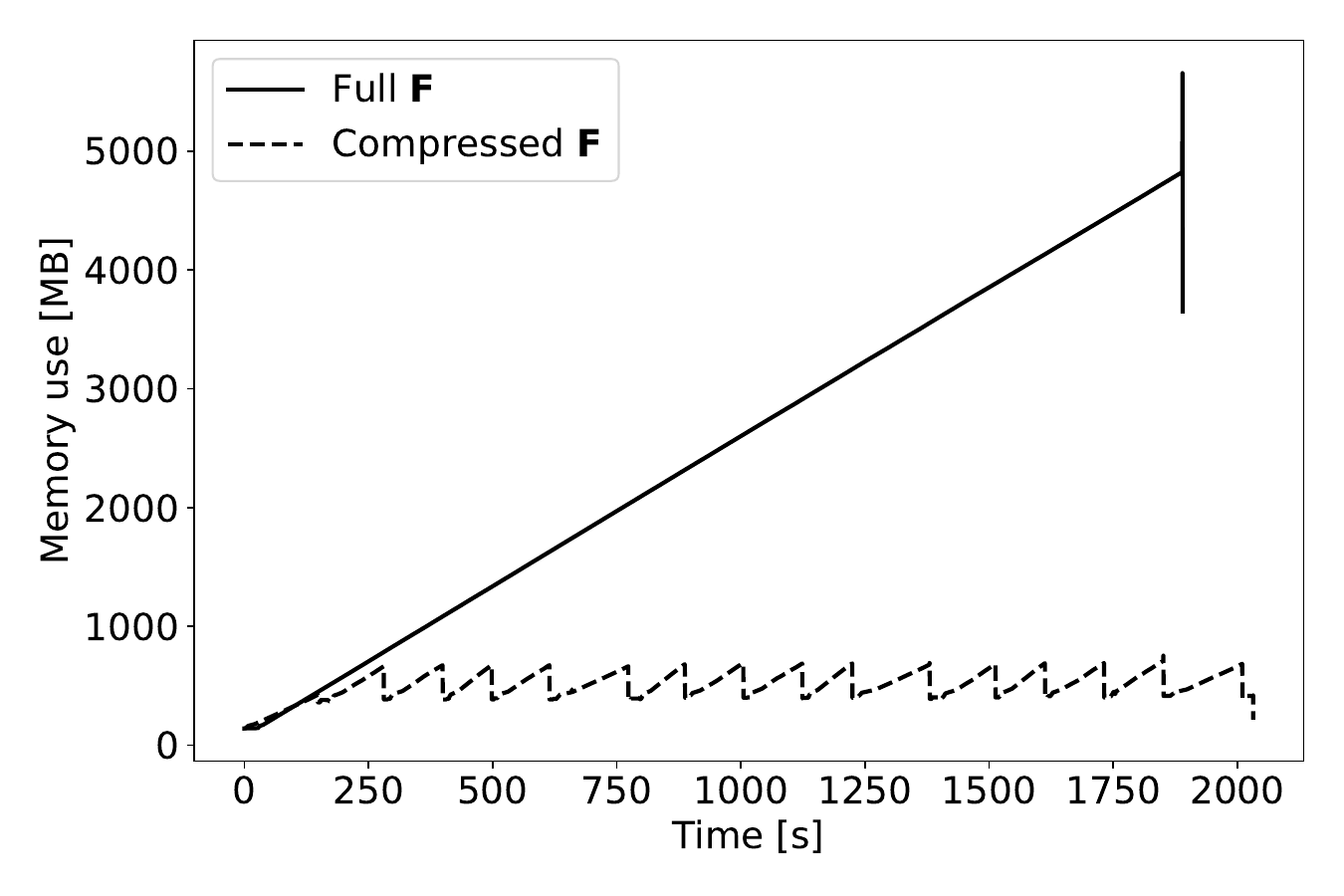}
  \caption{Memory use over time during view factor matrix assembly for $\m{F}$ and $\hat{\m{F}}$ for $N = 3.3 \times 10^4$ and $\varepsilon = 10^{-2}$. \emph{Note}: the sharp vertical ``up and down'' is not an error bar, and its source is uncertain. It could be an artifact of how we are measuring memory use, or it could be an artifact of how Python manages memory. We used a memory profiling package to collect the memory use statistics~\cite{MemoryProfiler}.}\label{fig:spherical-crater-memory-usage}
\end{figure}

To triangulate the hemispherical crater, we choose an edge length $h > 0$ and, depending on whether we discretize the shadow line, discretize the silhouette line and crater rim into three piecewise linear curves with segments of equal length: one discretizing $(x_s(y_s), y_s)$ with $-y_p \leq y_s \leq y_p$, and two more discretizing the crater rim, one on either side of $(x_s(y_s), \pm y_s)$. We do this so that the length of each segment in these curves is $O(h)$. Once we have these curves, we use Triangle~\cite{Shewchuk:1996aa} to mesh the horizontal plane, the shadowed part of the crater, and the illuminated part of the crater separately, concatenating the resulting meshes together. We uniformly refine the mesh so that the area of each 3D triangle is no greater than $\tfrac{2}{3}h^2$ (see Fig.\ \ref{fig:ingersoll-meshing}), which results in a mesh with $N$ faces.

For our numerical tests, we compare the compressed view factor matrix computed using Algorithm~\ref{algo:FF-assembly}, with compression tolerance $\varepsilon$ and $N$ triangles, which we denote $\hat{\m{F}}$, and the corresponding full view factor matrix computed using Algorithm~\ref{algo:FF-csr}, and denoted $\m{F}$. There are multiple sources of discretization error:
\begin{enumerate}
    \item The error due to the approximation of the true surface with a triangle mesh with planar facets.
    \item The error due to our point collocation approximation (using constant basis functions) of the integral equation used to compute $\m{F}$.
    \item The error due to approximating the view factors using the midpoint rule.
    \item The error due to approximating $\m{F}$ by hierarchically compressing it, resulting in $\hat{\m{F}}$.
\end{enumerate}
Our focus is on the fourth and final source of error.

For the spherical cap crater, we set $\esun = 15^\circ$, $\Ssun = 1000$ W/m$^2$, $\alpha = 0.3$, $\epsilon = 0.99$, $\beta = 40^\circ$, and $r_c = 0.8$ m. We construct a triangle mesh where the shadow line has been contoured exactly and one where it has not. To measure the error, we compute the relative $\ell_1, \ell_2$, and $\ell_\infty$ error of $\m{T}$ for varying $N$ and $\varepsilon$, as shown in \Cref{fig:spherical-crater-errors}---this figure demonstrates the convergence of the method. To evaluate the performance characteristics of our method, we time the assembly of $\m{F}$ and $\hat{\m{F}}$; i.e., we compare the time it takes to run the Algorithms~\ref{algo:FF-csr} and~\ref{algo:FF-assembly}. Next, we compare the amount of time it takes to compute the equilibrium temperature for $\m{F}$ and for each $\hat{\m{F}}$, as well as the size of each view factor matrix (see \Cref{fig:spherical-crater-performance-characteristics}). The time for the assembly of the compressed $\m{F}$ is approximately proportional to $O(N^2)$. The size of $\m{F}$ in memory is subproportional to $O(N^2)$. Using the compressed view factor matrix, the time to compute the terrain irradiance is approximately proportional to $O(N^2)$ instead of $O(N^3)$, which is the main result of this work. Finally, we measure the resident memory use during the construction of $\m{F}$ and $\hat{\m{F}}$ with $N = 3.3 \times 10^4$ and $\varepsilon = 10^{-2}$ in \Cref{fig:spherical-crater-memory-usage}. Our simulations were done using single-precision floating-point arithmetic.

For the range of $\varepsilon$ chosen, we identify three different qualitative behaviors (see \Cref{fig:spherical-crater-errors,fig:spherical-crater-performance-characteristics}). For $\varepsilon = 10^{-1}$ (the least accurate), the compressed view factor matrix is extremely small and fast, but does not approximate the full view factor matrix consistently, since it appears to have a degraded rate of convergence. For $\varepsilon = 10^{-2}$, the view factor matrix is about an order of magnitude larger, but it retains an MVP which appears to scale roughly as $O(N)$. Furthermore, it appears to approximate the full view factor matrix consistently, with only a modest amount of error incurred. However, for $\varepsilon = 10^{-3}$, we break the bank: too many resources are spent approximating the view factor matrix, and the linear time and space scalings break down.

\section{A realistic application in planetary science: de Gerlache crater}\label{sec:gerlache}

To test our numerical method on a more realistic problem, we run a time-dependent simulation of de Gerlache crater at the lunar South Pole ($88.5^\circ$S, $87.1^\circ$W). This crater hosts a large PSR which acts as a cold trap for volatiles~\citep{paige10a}. Small exposures of surface ice were observed from hyperspectral reflectance data~\citep{shuaili2018}. Parts of its rim are highly illuminated~\citep{mazarico11}, making it a potential site of interest for future robotic and human exploration. The de Gerlache crater has a complex morphology, of interest to test real-world applications of our algorithm. Part of its irregular floor is at a lower elevation which creates strong variations in the temperature field of the PSR.

To mesh de Gerlache, we use a recent NASA Digital Elevation Model (DEM)~\cite{barker2021improved} with a resolution of 5 m/pixel, and an elevation accuracy of 30-55 cm. The DEM itself was originally in a south polar stereographic projection. To make a triangle mesh, we select a circular region of interest (ROI) centered on de Gerlache, along with a larger, square bounding box. We set a background maximum triangle area of 3.0 km$^2$, and prescribe a smaller maximum area inside the circular ROI. We then use Shewchuk's Triangle~\cite{Shewchuk:1996aa} to create a Delaunay triangle mesh which satisfies these area constraints. This allows us to cheaply incorporate occlusion and self-heating by far-off topography directly into our view factor matrix.

\subsection{Numerical parameter study for de Gerlache crater}

\begin{figure}[tbh!]
 \begin{tabular}{cc}
 (a) & (b) \\
  \includegraphics[width=0.49\linewidth]{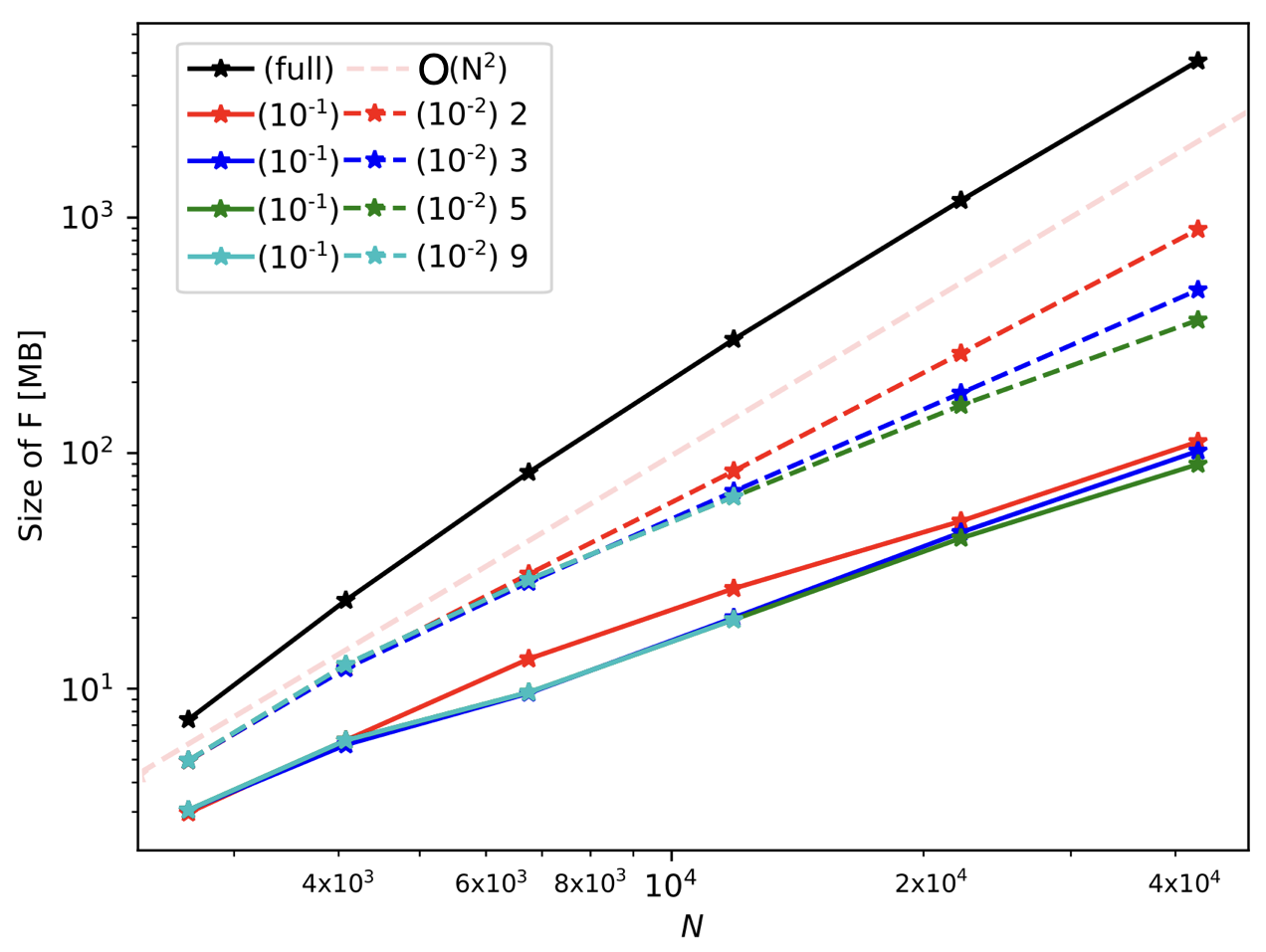}
  &
  \includegraphics[width=0.49\linewidth]{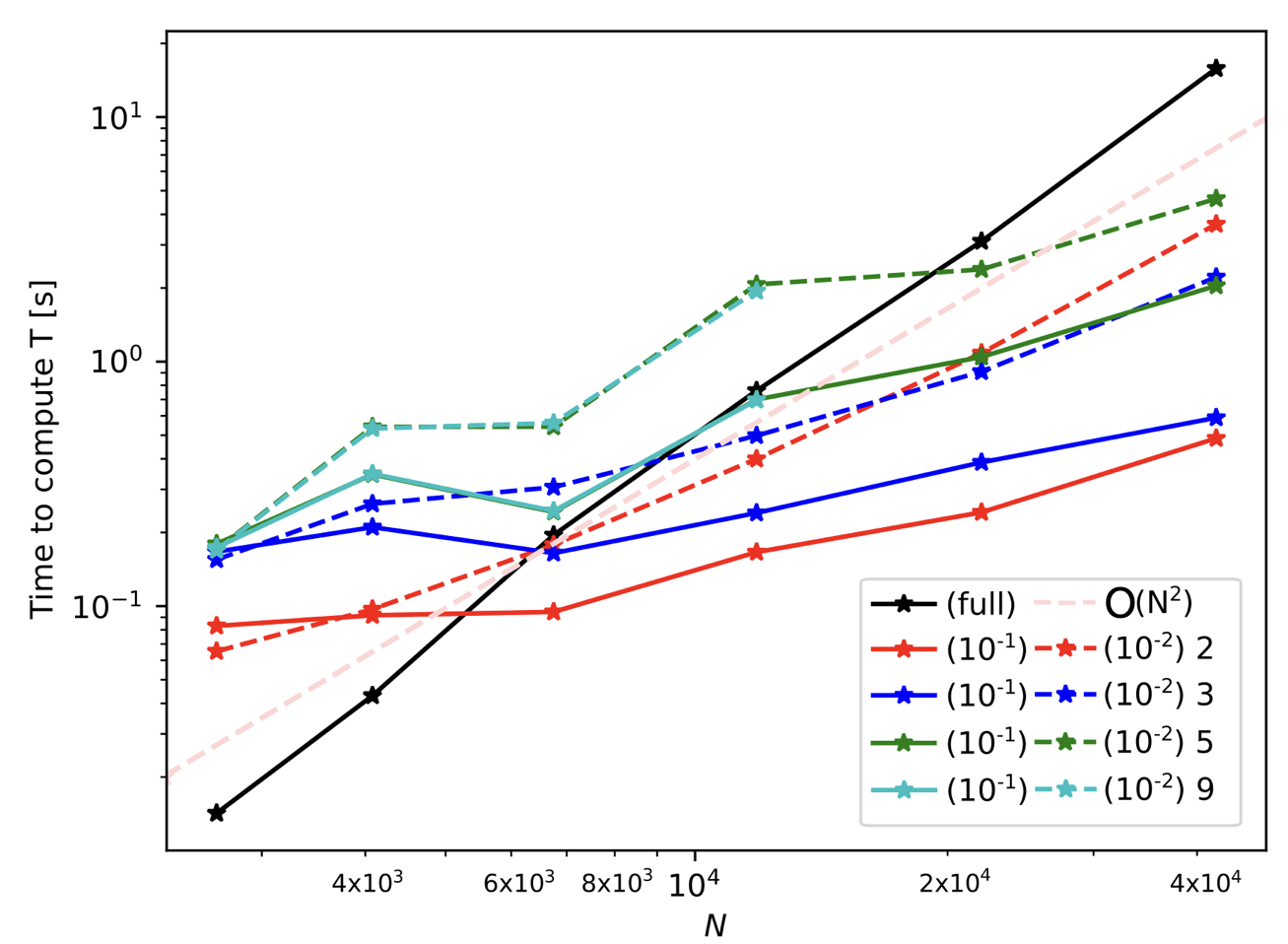}
  \end{tabular}
  \caption{Sizes and timings for thermal computations at the De Gerlache crater. The inner maximum area is set to $A = 1.6, 0.8, 0.4, 0.2, 0.1$, and $0.05$ km$^2$, while the outer maximum area is kept fixed at 3.0 km$^2$. The number of triangles $N$ is plotted on the $x$-axis, while the legend indicates ($\varepsilon$) and max depth for each curve. The hierarchically compressed view factor matrix is assembled using Algorithm~\ref{algo:FF-assembly} for varying $\varepsilon$ and maximum depths. \emph{Left}: the size of $\m{F}$ in MB. \emph{Right}: the timings of MVP operations necessary for computing surface temperature with $\m{F}$.}\label{fig:gerlache-performance-characteristics}
\end{figure}

\begin{figure}[tbh!]
  \includegraphics[width=\linewidth]{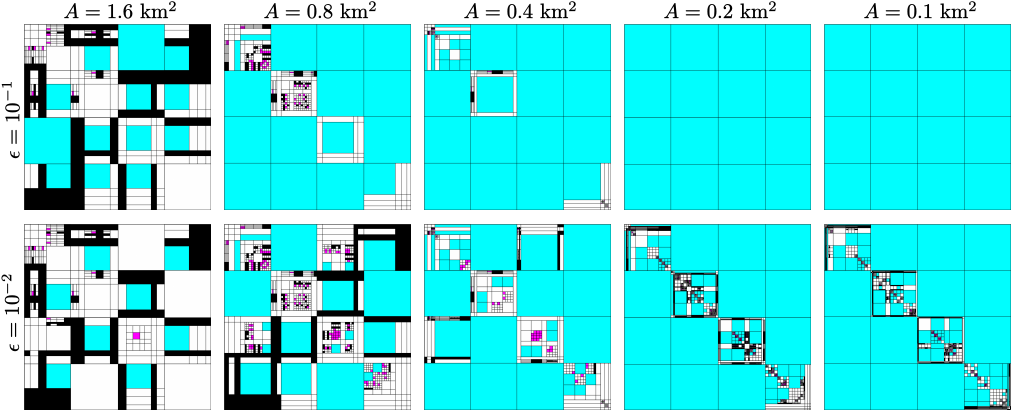}
  \caption{The block structure of the compressed view factor matrices for de Gerlache crater, for varying maximum triangle area inside the ROI ($A$) and compression tolerance ($\varepsilon$). The color scheme is as follows. {\color{cyan} \textbf{Cyan}}: sparse SVD block, {\color{black} \textbf{black}}: zero block, {\color{gray} \textbf{white}}: sparse CSR block, {\color{magenta} \textbf{magenta}}: dense block.}\label{fig:gerlache-blocks}
\end{figure}

Algorithm~\ref{algo:FF-assembly} hierarchically compresses $\m{F}$, recursively moving from one level of a spatial tree to the next until a user-defined block-size or maximum depth is reached. This implicitely defines the total number of leaf nodes. A large number of small blocks results in a smaller size and a reduced the number of FLOPs in the resultant MVP (see algorithm~\ref{algo:FF-multiply}). On the other hand, they also increase the assembly time and the complexity of the hierarchical matrix, which impacts the MVP efficiency in our implementation. This results in a fair amount of time spent passing data back and forth between Python and the lower-level libraries (such as BLAS) which implement the MVPs for individual blocks. In this section, we consider the effect of using a fixed level of refinement in the spatial tree to control the number of small blocks which are generated.

On an 8-core machine, we collect timing and memory measurements for the compressed view factor matrices for de Gerlache crater. We also collect timings for computing surface temperatures. Figure 5 shows that the time for each MVP is roughly proportional to $O(N)$ in a realistic application, compared to $O(N^2)$ if we use the uncompressed view factor matrix. We plot the blocks of $\m{F}$ for varying maximum inner areas ($A$) and compression tolerances ($\varepsilon$) in \Cref{fig:gerlache-blocks} to get a sense of the impact the compression algorithm has on the structure of the hierarchical matrix.

While the size of $\m{F}$ mainly depends on the chosen tolerance value, the speed of the MVP also strongly depends on the chosen maximum depth of the quadtree. The deeper the tree, the closer the MVP will come to scaling asymptotically like $O(N\log N)$, although the asymptotic prefactor will naturally be larger. The latter is a shortcoming of our current Python implementation, rather than an intrinsic limitation of our approach, and will be addressed in future works. When choosing an appropriate maximum tree depth for the problem at hand, our approach significantly outperforms the uncompressed form factor matrix. Since modeling scattered light and temperature on planetary surfaces involves a large number of MVPs (depending on the time span being analyzed and on the desired time resolution), order of magnitude gains for such simulations are potentially game-changing.

\subsection{Solving the time-dependent problem for de Gerlache crater}

\begin{figure}
  \includegraphics[width=\linewidth]{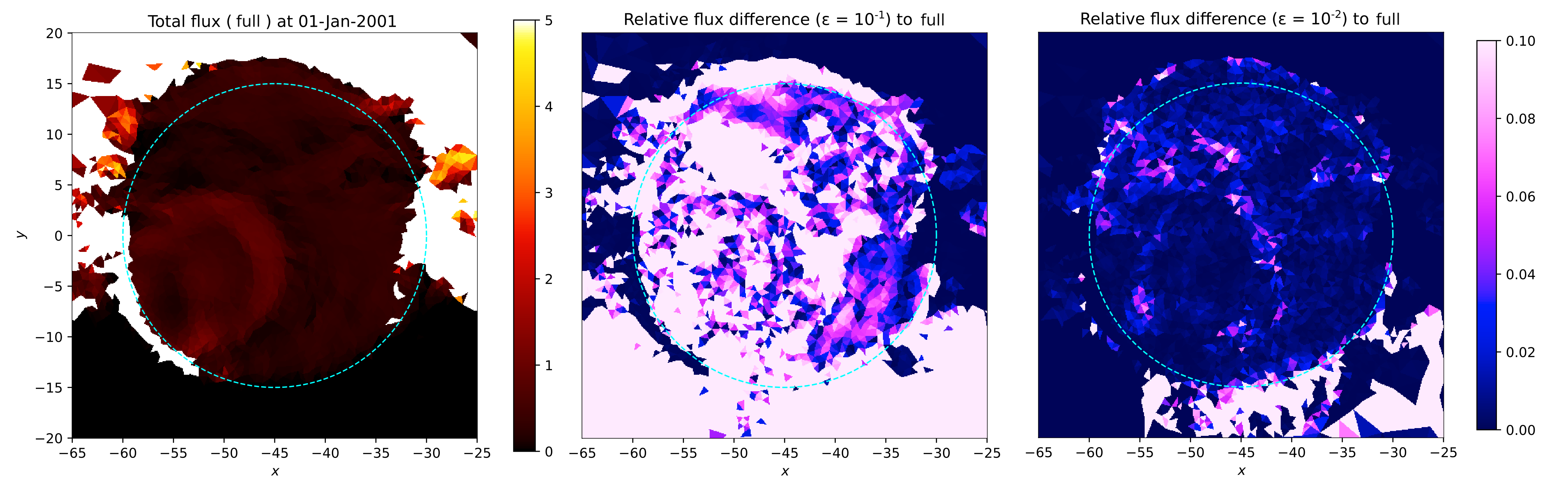}
  \caption{An overview of the total flux (direct and reflected, short and long wavelenghts) at de Gerlache crater at the Lunar south pole on 01-Jan-2001, with circular ROI imposed. We include a ``skirt'' of lower resolution triangles in the same watertight triangle mesh surrounding the circular ROI to incorporate occlusion and self-heating from far-off topography. For this mesh, we set the maximum triangle area outside the ROI to 3 km$^2$ and the maximum triangle area inside the ROI to $0.1$ km$^2$. \emph{Left}: the total flux computed using $\m{F}$ assembled with Algorithm~\ref{algo:FF-csr}, clamped to $5$ W/m$^2$ to highlight details of the crater floor. \emph{Middle ($\varepsilon = 10^{-1}$) and Right ($\varepsilon = 10^{-2}$)}: the pointwise absolute relative error $|\hat{\m{Q}} - \m{Q}|/|\m{Q}|$, with $\m{Q}$ computed from \eqref{eq:globalsolution} and $\hat{\m{Q}}$ computed using $\m{F}$ assembled with Algorithm~\ref{algo:FF-assembly}, but for varying $\varepsilon$.  As $\varepsilon$ is decreased, we see a concomitant reduction in the relative pointwise error, as expected. In each plot, the maximum relative error is clamped to 10\%, to show this effect.}\label{fig:gerlache-pointwise-Q}
\end{figure}

As an example, we follow Algorithm~\ref{algo:time-dependent} to perform a time-dependent simulation at the De Gerlache crater. We use the SPICE Toolkit through the SpiceyPy wrapper~\cite{annex2020spiceypy} to generate a set of Sun-crater geometries over several full illumination cycles. We start our simulation assuming a constant temperature of $110$ K for all mesh faces, an albedo $\alpha=0.2$, emissivity $\epsilon=0.95$, $\rho c = 9.6 \times 10^{-5}$ J/m$^3$K, and geothermal flux $F_g = 0.005$ W/m$^2$. We also consider $60$ subsurface layers of increasing thickness down to a depth of $2.5$~m (similar to what is done in~\cite{Schorghofer2021}), and let them reach thermal equilibrium (i.e., when inter-cycle temperature variations at all layers are below a threshold of $1$K). For this example, we end up dealing with $\sim 2 \times 10^4$ time steps when considering $20$ cycles and a time resolution of $0.3$ days. 
This is a computationally demanding thermal modeling of planetary surfaces, especially when dealing with high resolution datasets, and one where our compression algorithm shows clear advantages: indeed, the overhead of compressing $\m{F}$ only applies once, while more efficient MVPs improve the efficiency of the individual time steps and make the problem treatable on regular machines.
At each time-step we compute the direct flux from the Sun to the mesh and the reflected flux (including short and long wavelengths) between facets. We provide a snapshot in Fig.~\ref{fig:gerlache-pointwise-Q}, also showing the relative discrepancies due to matrix compression when using Algorithm~\ref{algo:FF-assembly}. These errors then propagate to the modeling of surface temperatures, so that tolerance $\varepsilon$ needs to be adapted in a trade-off between the required accuracy and computational speed.

Once equilibrium is reached, we evaluate the maximum and average surface temperature reached over a full cycle, which are key in assessing which molecules can survive at a given depth over the crater floor, as well as their state~\cite{Zhang2009,Schorghofer2021}.

\begin{figure}[!ht]
  \includegraphics[width=\linewidth]{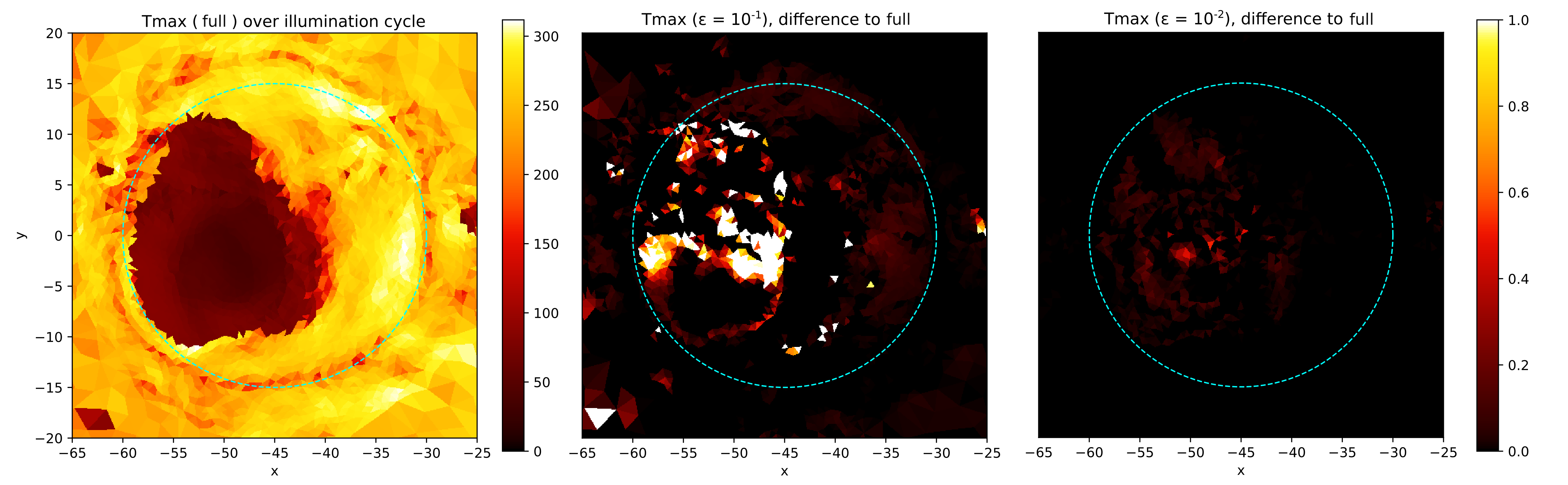} \\%
  \includegraphics[width=\linewidth]{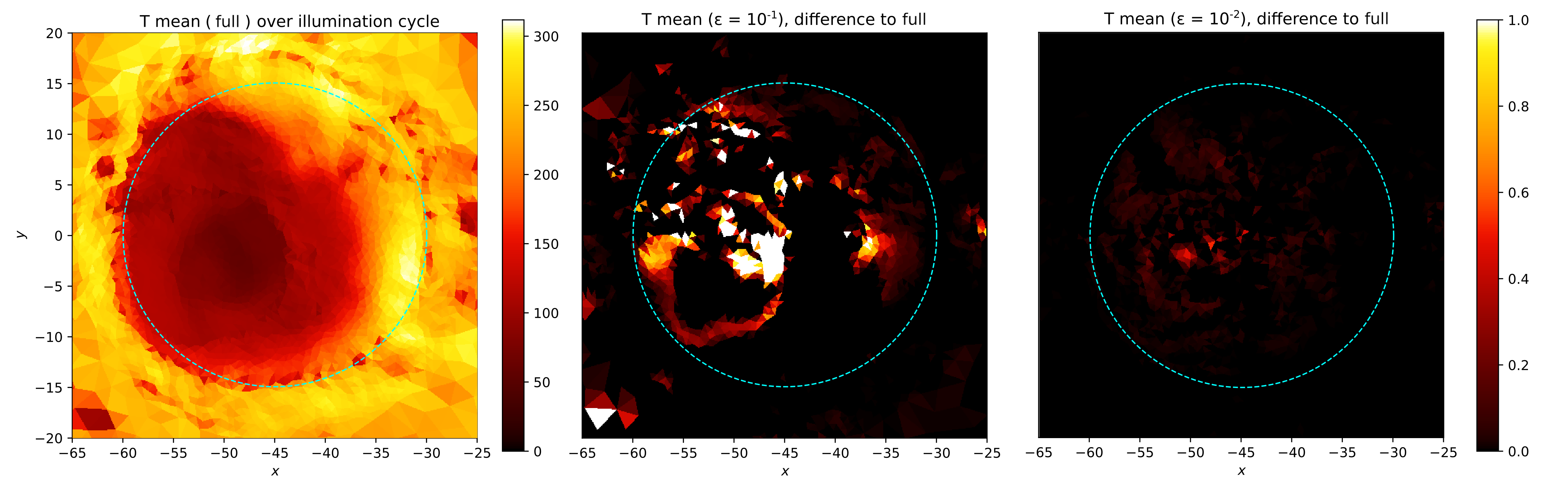}%
  \caption{Maximum and average temperature (K) at de Gerlache over a full illumination cycle (one year using a synthetic circular orbit of the Moon around the Sun) and at thermal equilibrium  and errors in the surface temperature
  for $\epsilon=10^{-1}$ and $\epsilon=10^{-2}$, computed using Algorithm~\ref{algo:time-dependent}. \emph{Top}: the maximum surface temperature, revealing the colder PSR. \emph{Bottom}: the average temperature over time, which is related with the sublimation rate of near-surface volatiles. All plots are in South Polar stereographic projection, km.}\label{fig:gerlache-year-stat-equil}
\end{figure}

Fig.~\ref{fig:gerlache-year-stat-equil} then shows the maximum and mean surface temperatures over a full illumination cycle as well as errors introduced by compressing $\m{F}$ with different tolerances.
We observe that while some drift between the ``true'' solution computed using the uncompressed view factor matrix and the approximate solution occurs, the bulk statistics agree well. In particular when using $\varepsilon=10^{-2}$, discrepancies are well below $1$K. Furthermore, the mesh fineness has a larger impact on the error than the compression tolerance.

\section{Global problems}

\begin{figure}[tbh!]
  \begin{subfigure}[b]{0.333\linewidth}
    \includegraphics[width=\textwidth]{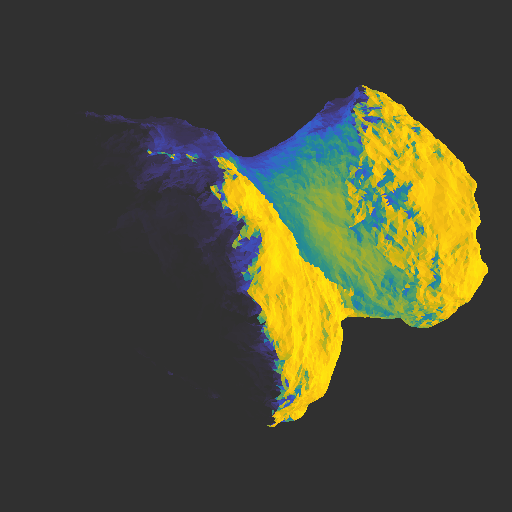}
    \caption{$\sim 49$K faces}
  \end{subfigure}
  \begin{subfigure}[b]{0.333\linewidth}
    \includegraphics[width=\textwidth]{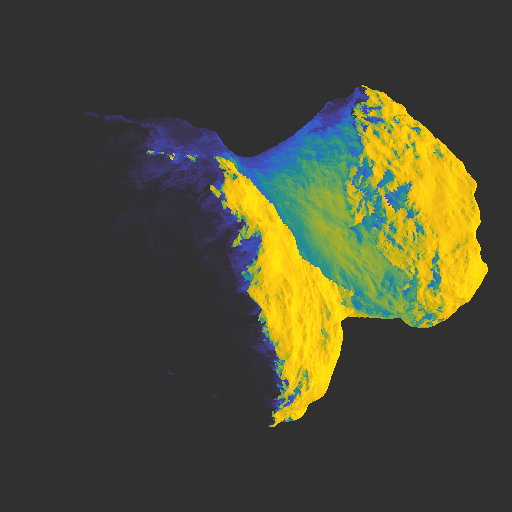}
    \caption{$\sim 98$K faces}
  \end{subfigure}
  \begin{subfigure}[b]{0.333\linewidth}
    \includegraphics[width=\textwidth]{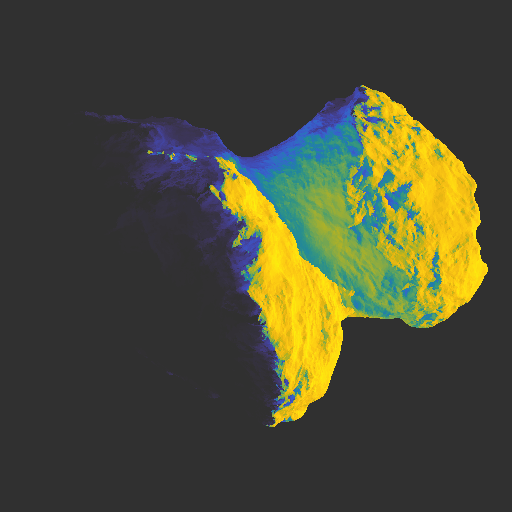}
    \caption{$\sim 197$K faces}
  \end{subfigure}
  \caption{Plots of the equilibrium temperature computed for one sun position for 67P/Churyumov-Gerasimenko, computed using the hierarchically compressed view factor matrices shown in \Cref{fig:67p-blocks}. As in \Cref{fig:67p-blocks}, fineness is increasing from left to right. In regions near the transition from sun to shadow, the temperature field resolves with the mesh. In regions receiving more direct sunlight, the temperature field is smoother and appears consistent, which agrees with \Cref{fig:gerlache-pointwise-Q}. Empirically, our method appears to be consistent.}\label{fig:67p-T}
\end{figure}

\begin{figure}[tbh!]
  \begin{subfigure}[b]{0.333\linewidth}
    \includegraphics[width=\textwidth]{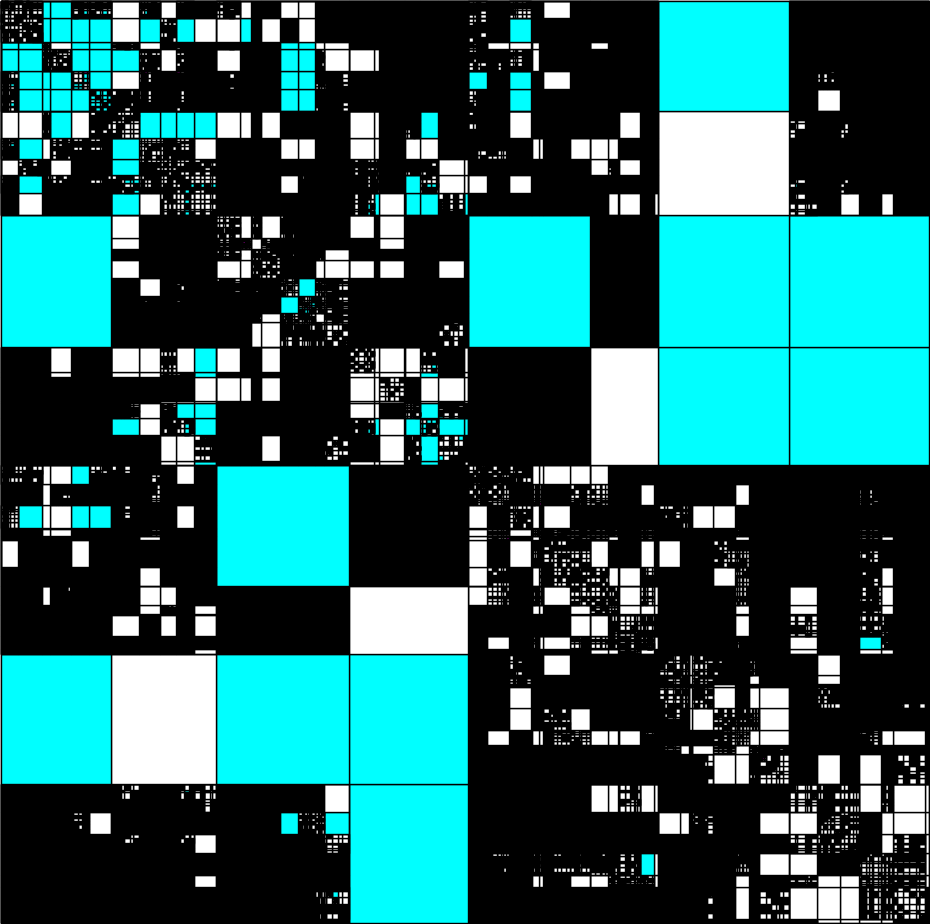}
    \caption{$\sim 49$K faces}
  \end{subfigure}
  \begin{subfigure}[b]{0.333\linewidth}
    \includegraphics[width=\textwidth]{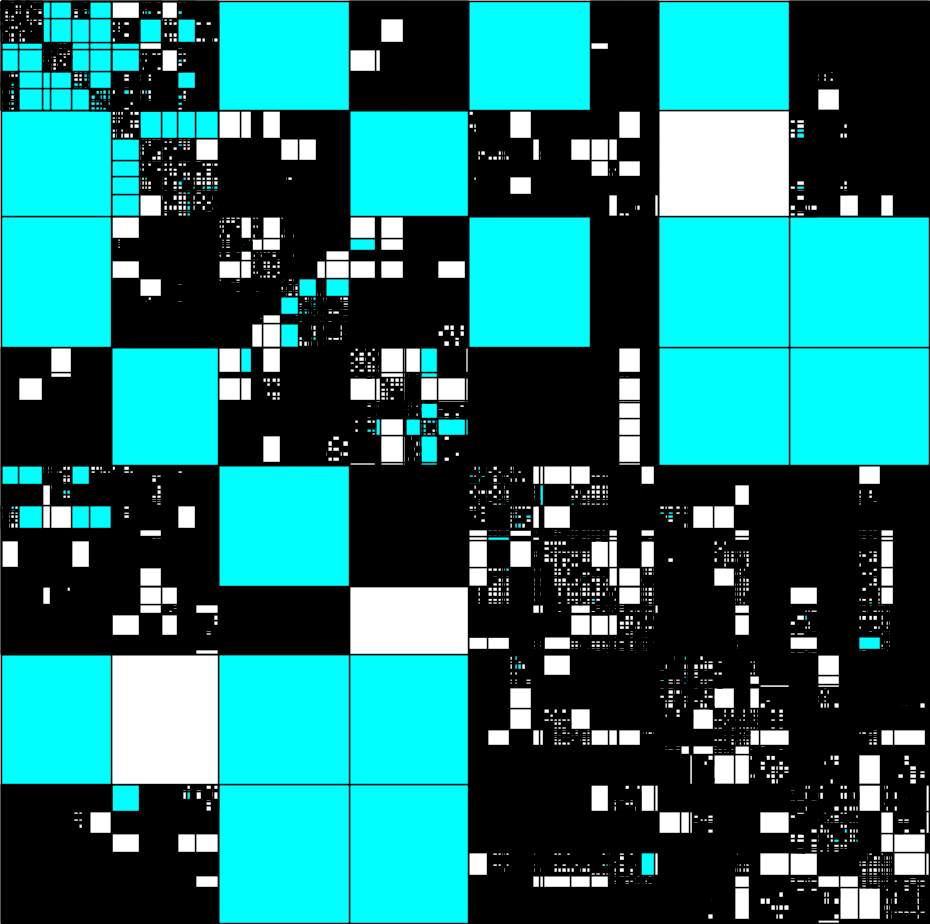}
    \caption{$\sim 98$K faces}
  \end{subfigure}
  \begin{subfigure}[b]{0.333\linewidth}
    \includegraphics[width=\textwidth]{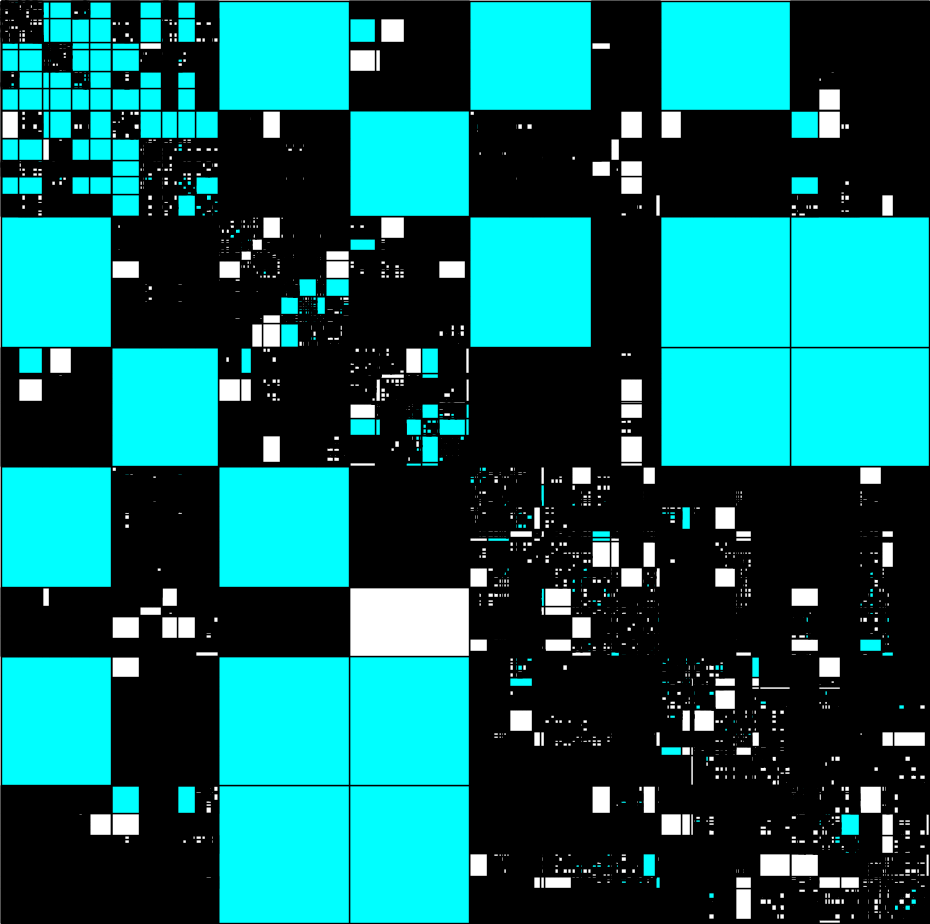}
    \caption{$\sim 197$K faces}
  \end{subfigure}
  \caption{Plots of the hierarchical view factor matrices $\m{F}$ for 67P/Churyumov-Gerasimenko for $\varepsilon = 10^{-2}$. See \Cref{fig:gerlache-blocks} for the color scheme. The fineness increases from left to right. As we increase the mesh fineness, we can see that CSR subblocks are gradually replaced by truncated SVD blocks stored in the CSR format.}\label{fig:67p-blocks}
\end{figure}

Finally, we run some tests applying our method to modeling radiation heat transfer on entire planetary bodies. We use a high-resolution shape model collected from Rosetta/OSIRIS observations~\cite{jorda2016global}. We downsample the mesh dyadically to obtain coarser test meshes of varying sizes using MeshLab's quadric edge decimation feature~\cite{garland1997surface,cignoni2008meshlab}. We use an octree to compute the hierarchically compressed view factor matrix $\m{F}$ for each of these meshes, and then compute the equilibrium temperature $\Teq$ for a single Sun position. See \Cref{fig:67p-T} to get a sense of the overall consistency of the method and the quality of the temperature field. Similar to \Cref{fig:gerlache-pointwise-Q}, the method appears to converge with the chosen tolerance of $\varepsilon = 10^{-2}$.

We also plot the sub-blocks of each $\m{F}$ matrix to see the manner in which the blocks refine. Compared to the $\m{F}$ matrices for the de Gerlache crater shown in \Cref{fig:gerlache-blocks}, the nature of the ``global'' problem reflected in \Cref{fig:67p-blocks} can be inferred by the higher percentage of zeros in $\m{F}$. Unlike the ``local''/``high-resolution'' problem, where nearly all entries are nonzero, here we get a large number of disconnected blocks, corresponding to pairs of octree cells with a relatively high degree of visibility.

\section{Conclusion and future work}

Large-scale thermal modeling on the surface of rough airless planetary bodies environment is a computationally demanding problem.
The number of element-to-element fluxes to calculate scales quadratically, and time-dependent calculations are often desired for numerous parameter combinations (scenarios).
Applying the radiosity method na\"{i}vely at each time step of a long-running simulation limits the size of problems which can be approached.

Borrowing ideas from the literature on hierarchical matrix factorizations and $n$-body codes, we present an algorithm for compressing the view factor matrix coming about as a result of discretizing the radiosity integral equation. It is not clear a priori whether this should be possible, due to the highly directional and sparse radiosity kernel. While the standard HODLR format fails to effectively compress the view factor matrix, a simple modification of it, which allows for off-diagonal blocks to be explored more deeply, succeeds. Unlike the HODLR format, this no longer guarantees $O(N \log N)$ time and space scaling---nevertheless, we observe close to linear scaling in our numerical experiments over a range of tolerances.

We demonstrate empirically that the compressed view factor matrix can be constructed in close to the amount of time it takes to assemble the original uncompressed version by using dynamic programming. Since the geometry (landscape) is static, once the view factor matrix has been assembled, it can be applied repeatedly. We find that for a reasonable set of tolerances, the time complexity scales nearly linearly in our numerical experiments. Additionally, only $O(N)$ space is required to store it in memory. In our numerical experiments, we find that this leads to orders of magnitude time and space savings. In these tests, we presented realistic examples showing how planetary scientists might use this method for modeling.

The algorithms here are presented as a prototype Python library. For speed, a future rewrite of the numerical kernel using C++ would be beneficial. The compressed view factor matrix is then the basis for (quicker) repeated computations of the scattered flux to get illumination and/or thermal conditions of the modeled surface at multiple epochs or under different thermophysical property assumptions.
Moreover, new missions and studies allow for ever larger, more detailed shape models of planetary surfaces. Such models are often too large for a typical workstation, and clusters have become a necessity for this kind of application. With our algorithm, one could first assemble---on a cluster---a significantly compressed version of a view factor matrix, but perform subsequent computations interactively on a laptop or workstation, thereby simplifying the model usability and workflow.

There are interesting directions which could be investigated as follow-on projects to the current work. Potentially, it would be possible to achieve a view factor MVP which is faster than $O(N^2)$ without an $O(N^2)$ precompute phase, which would result in concomitant improvements to the time complexity of Algorithm~\ref{algo:FF-assembly}. Developing fast direct solvers for our compressed view factor matrix is also of interest. These are algorithms in which a hierarchically compressed version of the inverse of a hierarchically compressed matrix is constructed. Although each solve already runs in $O(N)$ time due the radiosity system being extremely well-conditioned, a fast direct solver for this problem might allow one to carry out each time step even faster.

\section*{Acknowledgments}

This material is based upon work supported by the National Aeronautics and Space Administration through Grant Number 80NSSC19K0781 issued through the Internal Science Funding Model.

\bibliographystyle{elsarticle-num}
\bibliography{radiosity}

\begin{thebibliography}{10}
\expandafter\ifx\csname url\endcsname\relax
  \def\url#1{\texttt{#1}}\fi
\expandafter\ifx\csname urlprefix\endcsname\relax\def\urlprefix{URL }\fi
\expandafter\ifx\csname href\endcsname\relax
  \def\href#1#2{#2} \def\path#1{#1}\fi

\bibitem{martinsson2019fast}
P.-G. Martinsson, Fast Direct Solvers for Elliptic PDEs, SIAM, 2019.

\bibitem{arnold79}
J.~R. Arnold, Ice in the lunar polar regions, J. Geophys. Res. 84~(B10) (1979)
  5659--5668.
\newblock \href {https://doi.org/10.1029/JB084iB10p05659}
  {\path{doi:10.1029/JB084iB10p05659}}.

\bibitem{mazarico11}
E.~Mazarico, G.~A. Neumann, D.~E. Smith, M.~T. Zuber, M.~H. Torrence,
  Illumination conditions of the lunar polar regions using {LOLA} topography,
  Icarus 211 (2011) 1066--1081.

\bibitem{delbo15}
M.~Delbo, M.~Mueller, J.~P. Emery, B.~Rozitis, M.~T. Capria, Asteroid
  thermophysical modeling, in: P.~Michel, F.~{DeMeo}, W.~Bottke (Eds.),
  Asteroids IV, Univ. Arizona Press, Tucson, 2015, pp. 107--128.

\bibitem{mazarico18}
E.~Mazarico, M.~K. Barker, J.~B. Nicholas, {Advanced illumination modeling for
  data analysis and calibration. Application to the Moon}, Adv. Space Res.
  62~(11) (2018) 3214--3228.

\bibitem{rubanenko2018ice}
L.~Rubanenko, E.~Mazarico, G.~Neumann, D.~Paige, Ice in micro cold traps on
  mercury: Implications for age and origin, J. Geophys. Res. Planets 123~(8)
  (2018) 2178--2191.

\bibitem{glaser19}
P.~Gl{\"a}ser, D.~Gl{\"a}ser, Modeling near-surface temperatures of airless
  bodies with application to the {M}oon, Astronomy \& Astrophysics 627 (2019)
  A129.

\bibitem{schorghofer19}
N.~Schorghofer, J.~S. Levy, T.~A. Goudge, {High-resolution thermal environment
  of recurring slope lineae in Palikir Crater, Mars, and its implications for
  volatiles}, J. Geophys. Res. 124 (2019) 2852--2862.
\newblock \href {https://doi.org/10.1029/2019JE006083}
  {\path{doi:10.1029/2019JE006083}}.

\bibitem{king20}
O.~King, T.~Warren, N.~Bowles, E.~Sefton-Nash, R.~Fisackerly, R.~Trautner, {The
  Oxford 3D thermophysical model with application to PROSPECT/Luna 27 study
  landing sites}, Planet. Space Sci. 182 (2020) 104790.

\bibitem{hayne2021micro}
P.~O. Hayne, O.~Aharonson, N.~Sch{\"o}rghofer, Micro cold traps on the {M}oon,
  Nature Astronomy 5~(2) (2021) 169--175.

\bibitem{mahanti22}
P.~Mahanti, T.~J. Thompson, M.~S. Robinson, D.~C. Humm, View factor-based
  computation of secondary illumination within lunar permanently shadowed
  regions, IEEE Geoscience and Remote Sensing Letters 19 (2022) 1--4.

\bibitem{paige10a}
D.~A. Paige, et~al., Diviner {L}unar {R}adiometer observations of cold traps in
  the {M}oon's south polar region, Science 330 (2010) 479--482.

\bibitem{cohen1993radiosity}
M.~F. Cohen, J.~R. Wallace, P.~Hanrahan, Radiosity and realistic image
  synthesis, Morgan Kaufmann, 1993.

\bibitem{araki2020radiosity}
S.~J. Araki, Radiosity view factor model for sources with general distribution,
  J. Comput. Phys. 406 (2020) 109146.

\bibitem{hanrahan91}
P.~Hanrahan, D.~Salzman, L.~Aupperle, A rapid hierarchical radiosity algorithm,
  Computer Graphics 25 (1991) 197--206.

\bibitem{gortler1993wavelet}
S.~J. Gortler, P.~Schr{\"o}der, M.~F. Cohen, P.~Hanrahan, Wavelet radiosity,
  in: Proceedings of the 20th annual conference on Computer graphics and
  interactive techniques, 1993, pp. 221--230.

\bibitem{kahler2008h}
U.~K{\"a}hler, H\^2-wavelet {Galerkin} {BEM} and its application to the
  radiosity equation, Ph.D. thesis, Technische Universit\"at Chemnitz (2008).

\bibitem{atkinson2000numerical}
K.~Atkinson, D.~Chien, J.~Seol, Numerical analysis of the radiosity equation
  using the collocation method, Electronic Transactions on Numerical Analysis
  11 (2000) 94--120.

\bibitem{bekaert2001hierarchical}
P.~Bekaert, Hierarchical and stochastic algorithms for radiosity, Ph.D. thesis,
  Katholieke Universiteit Leuven (Belgium) (2001).

\bibitem{fernandez09}
E.~Fernández, Low-rank radiosity, in: F.~Ser{\'o}n, O.~Rodr{\'\i}guez,
  J.~Rodr{\'\i}guez, E.~Coto (Eds.), IV Iberoamerican Symposium in Computer
  Graphics SIACG, 2009.

\bibitem{aguerre2016hierarchical}
J.~P. Aguerre, E.~Fern{\'a}ndez, A hierarchical factorization method for
  efficient radiosity calculations, Computers \& Graphics 60 (2016) 46--54.

\bibitem{xia2010fast}
J.~Xia, S.~Chandrasekaran, M.~Gu, X.~S. Li, Fast algorithms for hierarchically
  semiseparable matrices, Numerical Linear Algebra with Applications 17~(6)
  (2010) 953--976.

\bibitem{hackbusch2000sparse}
W.~Hackbusch, B.~N. Khoromskij, A sparse {$\mathcal{H}$}-matrix arithmetic,
  Computing 64~(1) (2000) 21--47.

\bibitem{samet1984quadtree}
H.~Samet, The quadtree and related hierarchical data structures, ACM Computing
  Surveys (CSUR) 16~(2) (1984) 187--260.

\bibitem{fabri2009cgal}
A.~Fabri, S.~Pion, {CGAL:} the computational geometry algorithms library, in:
  Proceedings of the 17th ACM SIGSPATIAL International Conference on Advances
  in Geographic Information Systems, 2009, pp. 538--539.

\bibitem{wald2014embree}
I.~Wald, S.~Woop, C.~Benthin, G.~S. Johnson, M.~Ernst, Embree: a kernel
  framework for efficient cpu ray tracing, ACM Transactions on Graphics (TOG)
  33~(4) (2014) 1--8.

\bibitem{ying2004kernel}
L.~Ying, G.~Biros, D.~Zorin, A kernel-independent adaptive fast multipole
  algorithm in two and three dimensions, Journal of Computational Physics
  196~(2) (2004) 591--626.

\bibitem{lehoucq1998arpack}
R.~B. Lehoucq, D.~C. Sorensen, C.~Yang, ARPACK users' guide: solution of
  large-scale eigenvalue problems with implicitly restarted Arnoldi methods,
  SIAM, 1998.

\bibitem{zhao2005adaptive}
K.~Zhao, M.~N. Vouvakis, J.-F. Lee, The adaptive cross approximation algorithm
  for accelerated method of moments computations of emc problems, IEEE
  transactions on electromagnetic compatibility 47~(4) (2005) 763--773.

\bibitem{schorghofer2017planetary}
N.~Sch{\"o}rghofer,
  \href{https://github.com/nschorgh/Planetary-Code-Collection}{Planetary-code-collection:
  Thermal and ice evolution models for planetary surfaces v1.1.4} (2017).
\newblock \href {https://doi.org/10.5281/zenodo.1001854}
  {\path{doi:10.5281/zenodo.1001854}}.
\newline\urlprefix\url{https://github.com/nschorgh/Planetary-Code-Collection}

\bibitem{van2011numpy}
S.~Van Der~Walt, S.~C. Colbert, G.~Varoquaux, The numpy array: a structure for
  efficient numerical computation, Computing in Science \& Engineering 13~(2)
  (2011) 22--30.

\bibitem{virtanen2020scipy}
P.~Virtanen, R.~Gommers, T.~E. Oliphant, M.~Haberland, T.~Reddy, D.~Cournapeau,
  E.~Burovski, P.~Peterson, W.~Weckesser, J.~Bright, et~al., Scipy 1.0:
  fundamental algorithms for scientific computing in python, Nature Methods
  17~(3) (2020) 261--272.

\bibitem{Shewchuk:1996aa}
J.~R. Shewchuk, Triangle: Engineering a 2{D} quality mesh generator and
  {D}elaunay triangulator, in: Workshop on Applied Computational Geometry,
  Springer, 1996, pp. 203--222.

\bibitem{buhl68}
D.~Buhl, W.~J. Welch, D.~G. Rea, Reradiation and thermal emission from
  illuminated craters on the lunar surface, J. Geophys. Res. 73 (1968)
  5281--5295.

\bibitem{ingersoll92}
A.~P. Ingersoll, T.~Svitek, B.~C. Murray, {Stability of polar frosts in
  spherical bowl-shaped craters on the Moon, Mercury, and Mars}, Icarus 100~(1)
  (1992) 40--47.

\bibitem{schroder1993closed}
P.~Schr{\"o}der, P.~M. Hanrahan, A closed form expression for the form factor
  between two polygons, Princeton University, 1993.

\bibitem{walton2002calculation}
G.~N. Walton, Calculation of obstructed view factors by adaptive integration,
  Tech. rep. (2002).

\bibitem{narayanaswamy2015analytic}
A.~Narayanaswamy, An analytic expression for radiation view factor between two
  arbitrarily oriented planar polygons, Intern. J. Heat Mass Transf. 91 (2015)
  841--847.

\bibitem{MemoryProfiler}
memory-profiler - {PyPI}, \url{https://pypi.org/project/memory-profiler/},
  accessed: 2022-09-04.

\bibitem{shuaili2018}
S.~Li, P.~G. Lucey, R.~E. Milliken, P.~O. Hayne, E.~Fisher, J.-P. Williams,
  D.~M. Hurley, R.~C. Elphic, Direct evidence of surface exposed water ice in
  the lunar polar regions, Proc. Nat. Acad. Sci. 115~(36) (2018) 8907--8912.
\newblock \href {https://doi.org/10.1073/pnas.1802345115}
  {\path{doi:10.1073/pnas.1802345115}}.

\bibitem{barker2021improved}
M.~K. Barker, E.~Mazarico, G.~A. Neumann, D.~E. Smith, M.~T. Zuber, J.~W. Head,
  Improved lola elevation maps for south pole landing sites: Error estimates
  and their impact on illumination conditions, Planet. Space Sci. 203 (2021)
  105119.

\bibitem{annex2020spiceypy}
A.~M. Annex, B.~Pearson, B.~Seignovert, B.~T. Carcich, H.~Eichhorn, J.~A.
  Mapel, J.~L.~F. Von~Forstner, J.~McAuliffe, J.~D. Del~Rio, K.~L. Berry,
  et~al., Spiceypy: A pythonic wrapper for the spice toolkit, Journal of Open
  Source Software 5~(46) (2020) 2050.

\bibitem{Schorghofer2021}
N.~Schorghofer, J.-P. Williams, J.~Martinez-Camacho, D.~A. Paige, M.~A.
  Siegler, Carbon dioxide cold traps on the moon, Geophys. Res. Lett. 48~(20)
  (2021).
\newblock \href {https://doi.org/10.1029/2021gl095533}
  {\path{doi:10.1029/2021gl095533}}.

\bibitem{Zhang2009}
J.~A. Zhang, D.~A. Paige, Cold-trapped organic compounds at the poles of the
  moon and mercury: Implications for origins, Geophys. Res. Lett. 36~(16)
  (2009).
\newblock \href {https://doi.org/10.1029/2009gl038614}
  {\path{doi:10.1029/2009gl038614}}.

\bibitem{jorda2016global}
L.~Jorda, R.~Gaskell, C.~Capanna, S.~Hviid, P.~Lamy, J.~{\v{D}}urech, G.~Faury,
  O.~Groussin, P.~Guti{\'e}rrez, C.~Jackman, et~al., The global shape, density
  and rotation of comet 67p/churyumov-gerasimenko from preperihelion
  rosetta/osiris observations, Icarus 277 (2016) 257--278.

\bibitem{garland1997surface}
M.~Garland, P.~S. Heckbert, Surface simplification using quadric error metrics,
  in: Proceedings of the 24th annual conference on Computer graphics and
  interactive techniques, 1997, pp. 209--216.

\bibitem{cignoni2008meshlab}
P.~Cignoni, M.~Callieri, M.~Corsini, M.~Dellepiane, F.~Ganovelli, G.~Ranzuglia,
  et~al., Meshlab: an open-source mesh processing tool., in: Eurographics
  Italian chapter conference, Vol. 2008, Salerno, Italy, 2008, pp. 129--136.

\end{thebibliography}

\end{document}